\documentclass[11pt,reqno]{amsart}

\usepackage{amsmath, amsthm, amssymb,amsaddr}
\usepackage{dcolumn}
\usepackage{bm}
\usepackage{amsmath} 
\usepackage{mathrsfs}
\usepackage{mathtools}
\usepackage{algpseudocode}
\usepackage{algorithm}
\usepackage[pdftex]{graphicx}
\usepackage{subfigure}
\usepackage{epstopdf}
\usepackage{tikz}
\usepackage{xcolor}
\usepackage[export]{adjustbox}
\usepackage{caption}
\usepackage{hyperref}

\theoremstyle{remark}

\newtheorem{definition}{Definition}

\DeclareMathOperator*{\argmin}{arg\,min}

\newcommand{\be}[1]{\begin{equation}\label{#1}}
\newcommand{\ee}{\end{equation}}

\newcommand{\no}[1]{#1}

\newcommand{\bx}{\boldsymbol{x}}

\newcommand{\bk}{\boldsymbol{k}}

\newcommand{\bn}{\boldsymbol{n}}
\newcommand{\iu}{{i\mkern1mu}}
\newcommand{\bkappa}{\boldsymbol{\kappa}}
\newcommand{\Lim}[1]{\raisebox{0.5ex}{\scalebox{0.8}{$\displaystyle \lim_{#1}\;$}}}

\renewcommand{\no}[1]{} 

\title{Ray-based inversion accounting for scattering for biomedical ultrasound tomography}

\author{Ashkan Javaherian and Ben Cox}
\email{a.javaherian@ucl.ac.uk}
\address{Department of Medical Physics \& Biomedical Engineering,\\ University College London, London, UK. WC1E 6BT}
\date{May 2021}

\begin{document}

\maketitle


\begin{abstract}
An efficient and accurate image reconstruction algorithm for ultrasound tomography (UST) is described and demonstrated, which can recover accurate sound speed distribution from acoustic time series measurements made in soft tissue. The approach is based on a second-order iterative minimisation of the difference between the measurements and a model based on a ray-approximation to the heterogeneous Green's function. It overcomes the computational burden of full-wave solvers while avoiding the drawbacks of time-of-flight methods. Through the use of a second-order iterative minimisation scheme, applied stepwise from low to high frequencies, the effects of scattering are incorporated into the inversion.
\end{abstract}

\section{Introduction}
The aim of ultrasound tomography (UST) is to estimate the acoustic properties of the interior of an object from ultrasonic measurements made on its boundary \cite{Ruiter,Hopp1}. There are two related steps: designing and constructing the measurement hardware, and designing and implementing the algorithms to reconstruct the images from the measured data \cite{Gemmeke}. This paper is concerned with the latter challenge, image reconstruction. Approaches to UST reconstruction can be classified by (1) whether the inversion is linearised or nonlinear, (2) the nature of the forward model used, and (3) the data type used in the inversion. In terms of the data used, there are broadly two categories. The first uses the complete measured time series, including the scattered waves, or its frequency domain components. The second uses only the direct \textit{times-of-flight} between the emitters and receivers; in other words, the data consists just of the time of the first arrival, and no scattered waves are included \cite{Duric,Li3,Li1,Li2,Op,Ali,Javaherian}. The approach we propose here for UST performs a nonlinear inversion by minimising the discrepancy between acoustic pressure time series data and a forward model based on a ray-approximation to the Green's function.

Recently, nonlinear inversion schemes based on a minimisation of the norm of the discrepancy between the measured time series and a full-wave numerical model, which depends on the unknown sound speed, have become popular. Because of the flexibility in the choice of model, this approach - which has become known as \textit{full-wave inversion} \cite{Plessix} - can be very general. When the degree of model-mismatch is low, ie.\ when the numerical model accurately represents the measurement scenario, these approaches have the potential to reconstruct accurate, high resolution images. This approach, for which there is a considerable literature in the seismic community \cite{Tromp}, has begun to be explored in earnest for medical applications \cite{Wiskin1,Wang,Goncharsky,Matthews,Matthews1, Wiskin2,Liva,Bachmann,Guasch2,Lucka}. Full-wave approaches are discussed further in Sec.\ \ref{sec:discussion}, but the biggest challenge with such schemes is the computational cost.

In this paper, we make contributions to both the forward and the inverse problems of UST. We avoid the computational challenge that full-wave models present by using a solver based on a ray-approximation to the Green's function. The ray-based forward model we present is a frequency domain model that can account for refraction by using bent rays, geometric spreading through Green's law, and arbitrary absorption (we use Szabo's model \cite{Szabo,Liebler,Kelly} but the formulation is general). The principle drawback of using a forward model based on rays is that it inherently neglects scattering. Overcoming this limitation is this paper's principle contribution to tackling the inverse problem. At each step, the Gauss-Newton search direction is computed by an iterative and implicit computation of the Hessian matrix. Solving the nonlinear inversion in this way implicitly accounts for the primary scattered field (which is sufficient for soft tissue) and can therefore provides high spatial resolution. The non-linearity of the inverse problem and the cycle-skipping problem are handled by first reconstructing a low-contrast sound speed image using a time-of-flight approach as an initial guess, and then solving the inverse problem from low to high frequencies \cite{Wiskin1,Wiskin2}. This combination of a ray-based forward model and a second-order inversion scheme provide a computationally efficient method for waveform tomography that avoids the principle drawbacks of both time-of-flight approaches and full-wave solvers.

In Sec.\ \ref{sec:forward-inverse-operators}, the forward and inverse problems of UST are introduced. In Sec.\ \ref{sec:greens} the forward model - the ray-based approximation to the heterogeneous Green's function - is described. Sec.\ \ref{sec:optimisation-greens} explains the second-order approach to the inverse problem, and the calculation of the gradient and Hessian based on the Green's function. The numerical tracing and discretisation of the rays and the Green's functions is described in Secs.\ \ref{sec:raytracing} and \ref{sec:Discretised-approximation}. Numerical examples demonstrating the methods in 2D are given in Sec.\ \ref{sec:num-results}, although all the results are applicable to 3D. 
A discussion of the significance of the results and links to similar work follows in Sec.\ \ref{sec:discussion}, and the paper ends with a brief conclusion section.

\section{Ultrasonic Sound Speed Tomography}
\label{sec:forward-inverse-operators}
This section describes the forward and inverse problems of sound speed tomography. Let $\bx= \left( x^1,...,x^d \right)$ denote a spatial position in $\mathbb{R}^d$ with $d$ the dimension. In general, $d$ can be either 2 or 3. This study is restricted to $d=2$. Accordingly, $\Omega \subset \mathbb{R}^d $ is an open bounded set, and contains the spatially-varying part of the sound speed distribution, $c(\bx)$, i.e. $(c(\bx)/c_0 - 1) \in C_0^\infty$, where $c_0$ is a scalar value representing the sound speed outside $\Omega$ (here the sound speed in water). 
Also, $\rho(\bx)$ will represent the spatially varying density. 
The open set $\Omega$ is bounded by a circular ring $\mathbb{S} \subset \mathbb{R}$ containing the emission and reception elements. 

\subsection{Forward problem}

\subsubsection{Excitation} 
Sequentially, each emission element, referred to here as emitter $e \in \left\{ 1,...N_e \right\}$, is excited by a pulse and acts as a source $s(t;\bx_e)$ within the excitation time $ t \in \big(0, T_s\big) $. The stack-vector of these source time series for all emission elements is denoted by $S \in  \mathbb{R}^{(0,T_s) \times N_e} $. Each emitter $e$ is idealised as a point source at $\bx_e$ with a directional dependence such that the induced acoustic pressure field, $p(t, \bx; \bx_e) \in  (0,\infty) \times \mathbb{R}^d $, will be a function of $(\bx - \bx_e)\cdot \bn_e$, where $\bn_e$ is a unit vector giving the orientation of the emitter. (As an example of a practical realisation, each emitter $e$ might be a finite-sized disc centered at $\bx_e$ with normal $\bn_e$.)

\subsubsection{Measurement} 
The induced acoustic pressure field is measured at the reception elements, referred to here as receiver $r \in \left\{ 1,...,N_r \right\} $, for times $t \in (0,T)$ with $T \gg T_s$ . For each excitation element $e$, the time series measured by the receiver $r$, centered at $\bx_r$, is represented by $ p{( t, \bx_r; \bx_e )} $, where
\begin{align}
\label{eq:measurement-operator}
\begin{split}
&\mathcal{M}_{(e,r)}: (0,\infty) \times\mathbb{R}^d  \rightarrow \mathbb{R}^{N_t}   \\
&p{( t, \bx_r; \bx_e )} = \mathcal{M}_{(e,r)}\big[p( t, \bx; \bx_e)\big],
\end{split}
\end{align}
where $N_t$ is the number of measurement time samples. We have used $t$ here to denote both the continuous and discrete time variables, but the meaning in any particular case should be clear from the context. Furthermore, what is actually measured at $\bx_r$ is an electrical signal, but we used the same notation $p$ for brevity. As well as acting as a sampling or discretisation operator, $\mathcal{M}_{(e,r)}$ incorporates a filtering stage, corresponding to the frequency-angle dependent response of the receiver. In other words, in general, $\mathcal{M}_{(e,r)}$ not only applies a frequency filter but also depends on $(\bx - \bx_r)\cdot \bn_r$, where $\bn_r$ is a unit vector giving the orientation of the receiver. (As with the source elements, a practical realisation of a receiver might be a finite-sized disc centered at $\bx_r$ with normal $\bn_r$.) 
While the data is typically measured in the time domain, ie. with a broadband excitation signal, the image reconstruction will be performed in the frequency domain. To this end, we define the following Fourier transform pair between the time and temporal frequency domains,
\begin{align}
p(\omega) = \mathscr{F}p(t) = \int_{-\infty}^{\infty} p(t) e^{i\omega t} dt,\quad
p(t) = \mathscr{F}^{-1}p(\omega) = \frac{1}{2\pi}\int_{-\infty}^{\infty} p(\omega) e^{-i\omega t} dt.
\end{align}

\begin{definition}
Here, the forward operator is defined in the frequency domain.
\begin{align}
\begin{split}
    \mathcal{A}:\mathbb{D}(\Omega) \rightarrow \mathbb{R}^{N_{\omega} N_r N_e}\\
    P = \mathcal{A} \big[c (\bx)\big],
\end{split}
\end{align}
where $N_{\omega}$ is the number of discrete frequencies used. Also, the space $\mathbb{D}$ is defined such that any function $c(\bx)\in \mathbb{D}$ satisfies $\big( c(\bx)/c_0 - 1 \big) \in C_0^\infty (\Omega) $, and $P \in \mathbb{R}^{N_{\omega} N_r N_s}$ is a stack-vector of measured complex amplitudes with components
\begin{align}  \label{eq:forward-operator}
    p_{( c; \omega, r; e )} = \mathcal A_{(\omega, r; e)} \big[c \big],
\end{align}
where $\mathcal A_{(\omega, r; e)}$ accounts for the Fourier transformation $\mathscr{F}$, acoustic propagation, the filtering and sampling described by $\mathcal{M}_{(e,r)}$ in the frequency domain, and the selection of the $N_{\omega}$ frequencies.
\end{definition}

\subsection{Inverse problem}
\label{sec:inverse_problem}
The inverse problem is now an estimation of the sound speed distribution $c(\bx)$ from the measured data $\hat{p}_{(\omega, r; e)}$, which is the set of recorded complex amplitudes for all pairs of emission and reception points, given the excitation amplitudes $s(\omega, \bx_e) = \mathscr{F}s(t,\bx_e)$ \cite{Plessix}. (Note that we use the same notations for parameters in time and frequency domains for brevity.) 
This can be posed as the nonlinear minimisation problem
\begin{align}   \label{minimisation-problem}
c^* = \argmin_{c}  \mathcal{F}(c),
\end{align}
where the objective function is in the form
\begin{align}   \label{objective-function}
     \mathcal{F}(c) = \frac{1}{2} \int \ \ \Big(  \sum_{e,r}  \left\lvert  \mathcal{A}_{(\omega, r; e)} \big[c \big] - \hat{p}_{(\omega, r;e)} \right\rvert ^2 \Big) \ d \omega .
\end{align}
The sound speed is inversely proportional to the wavenumber (see the next section), the refractive index, and the slowness. This inversion is therefore sometimes posed in terms of these alternative quantities.

\section{Approximate Green's Function for the Lossy Helmholtz Equation} \label{sec:greens}

This section describes how the propagation of acoustic waves are modelled in this study, and introduces the approximate Green's function which is the basis of the inversion method using in Sec. \ref{sec:optimisation-greens}.

\subsection{Lossy Helmholtz equation and complex wavevector}
A single frequency acoustic pressure field, $p(\omega,\bx)$, in an absorbing medium is often modelled using a lossy Helmholtz equation of the form
\begin{align} 
\label{eq:second-order-lossy-omega}
\left( \tilde{k}(\bx)^2 + \nabla^2  \right) p(\omega, \bx; \bx_e) =  -s(\omega,\bx_e), \quad \tilde{k} = |\tilde{\bk}|,
\end{align}
where $\tilde{\bk}$ is a complex wavevector 
that depends on the sound speed $c(\bx)$, the absorption coefficient $\alpha(\bx)$ and the frequency $\omega$. In our case $s$ represents a point source at position $\bx_e$. By considering a plane wave solution of the form
\begin{align}
\label{eq:plane_wave}
p(\omega, \bx)\propto \exp(i(\tilde{\bk}\cdot\bx - \omega t)),    
\end{align}
and writing $\tilde{\bk} = \bk + i \bk_i$, it is clear that the real part $\bk$ is related to the phase speed $c_p(\omega)$ of the wave by 
\begin{align} 
|\bk| \equiv k = \omega / c_p(\omega),
\end{align}
and the imaginary part is related to the absorption coefficient $\alpha$ by
\begin{align} 
\bk_i = \alpha (\bk/k).
\end{align}

\subsection{Dispersion relation}
The wavenumber $\tilde{k}$ and temporal frequency $\omega$ are related via a \textit{dispersion relation}, which encapsulates the effects of absorption and dispersion on the wave. In general, then
\begin{align} 
\label{eq:general_dispersion_relation}
k(\bx) = \zeta(\omega, c(\bx),\alpha(\bx)).
\end{align}
Many different models of absorption (and the accompanying dispersion) have been proposed for describing the behaviour of soft biological tissue, and each will lead to a different dispersion relation, a different form for the function $\zeta$. As an example, in this subsection we describe the particular absorption model used in the numerical examples in Sec. \ref{sec:num-results}, and the dispersion relation that results from it. 
However, it is important to note that the methods described in this section, Sec.\ \ref{sec:greens}, apart from this one subsection, do not depend on the particular form of the dispersion relation. 

\subsubsection{Szabo absorption model}
\label{sec:szabo_model}
A popular model for describing the propagation of acoustic waves in soft tissues is Szabo's absorbing wave equation \cite{Szabo}
\begin{align}  \label{eq:second-order-lossy}
       \left( \frac{1}{c(\bx)^2}  \frac{\partial^2 }{\partial t^2}  - \nabla^2  + \frac{2 \alpha_0(\bx) }{ c(\bx) \cos(\pi y/2)} \frac{\partial^{y+1}}{\partial t^{y+1}}  \right) p(t, \bx; \bx_e ) = s(t, \bx_e),
\end{align}
where $\partial^{y+1}/\partial t^{y+1}$ is a fractional time derivative resulting in an absorption coefficient $\alpha$ which follows the frequency power law $\alpha = \alpha_0 \omega^y$. Here, $\alpha_0$ has units $\text{Np} (\text{rad}/\text{s})^{-y}\text{m}^{-1}$, and $y$ is the power-law exponent with a non-integer often in the range $1 \leqslant y \leqslant 1.5$ for soft tissue \cite{Szabo, Liebler}. In the frequency domain, \eqref{eq:second-order-lossy} becomes the lossy Helmholtz equation, \eqref{eq:second-order-lossy-omega}, with complex wavenumber given by \cite{Kelly}
\begin{align}
\label{eq:szabo_dispersion_equation}
    \tilde{k} = k + i\alpha = \frac{\omega}{c} + \alpha \big( \tan{(\pi y/2)} +\iu  \big).
\end{align}

\subsection{Green's functions}
\label{sec:greens_functions}
In general, the solution of \eqref{eq:second-order-lossy-omega} can be written in terms of a Green's function as
\begin{align}   \label{eq:general-greens-source}
    p(\omega,\bx) = \int g(\omega, \bx; \bx') s(\omega,\bx') d\bx'
\end{align}
where the Green's function $g$ satisfies
\begin{align}  \label{eq:greens-wave-equation}
\left( \tilde{k}(\bx)^2 + \nabla^2  \right) g(\omega, \bx; \bx') =  -\delta(\bx - \bx').
\end{align}
As $s$ is a point source, this becomes simply a scaling
\begin{align}   \label{eq:pressure-greens-source}
    p(\omega,\bx; \bx_e) = g(\omega, \bx; \bx_e) s(\omega,\bx_e).
\end{align}

\subsubsection{2D case}
In the homogeneous, lossless case the two-dimensional (2D) free-space Green's function is given by \cite{Abramowitz, Barnett}
\begin{align}
\label{eq:greens-homogeneous}
    g_{0,2D}(\omega,\bx; \bx') &= \frac{\iu}{4} H_0^{(1)}  \Big(  \phi_0(\bx;\bx') \Big)\\
    & \approx A_0(\bx;\bx') \exp{\big( i(\phi_0(\bx;\bx') + \pi/4)  \big) },
\end{align}
where $A_0$ is an amplitude factor and $\phi_0$ is the phase. (The second line gives the large-argument limit approximation, for which the second terms and higher of an associated asymptotic expansion has been neglected, see \cite{Barnett} for further details.) In this lossless, homogeneous case, the phase speed $c_p = c_0$, the wavenumber $\tilde{k} = k_0 = \omega /c_0$ is real, and the phase is given by 
\begin{align} \label{eq:phase-homogenous}
  \phi_0(\bx;\bx') =  k_0 |\bx - \bx'|.
\end{align}
The amplitude is governed purely by geometric spreading (cylindrical spreading):
\begin{align} \label{eq:amplitude-homogenous}
A_0(\bx;\bx')
    &= \Big( 8 \pi k_0 |\bx - \bx'| \Big)^{-1/2} = \Big( 8 \pi \phi_0(\bx;\bx') \Big)^{-1/2}.
\end{align}
In the heterogeneous, absorbing, case we need an approximate form for the Green's function, so we write, by analogy with the homogeneous case:
\begin{align}
\label{eq:2D_greens}
    g_{2D}(\omega,\bx; \bx') & \approx A(\bx;\bx') \exp{\big( i(\phi(\bx;\bx') + \pi/4)  \big) },
\end{align}
where the amplitude factor now contains contributions from absorption and refraction as well as geometric spreading. 

\subsubsection{3D case}
In the 3D case the homogeneous, lossless Green's function is
\begin{align}
    g_{0,3D}(\omega,\bx; \bx') &= A_0(\bx;\bx') \exp{ \Big( i \phi_0(\bx;\bx') \Big)},
\end{align}
where $\phi_0 = k_0|\bx - \bx'|$ as for the 2D case but there is spherical rather than cylindrical spreading:
\begin{align}
A_0(\bx;\bx') = (4 \pi |\bx-\bx'| )^{-1} = k_0 ( 4 \pi \phi_0(\bx,\bx') )^{-1}.
\end{align}
So the approximate Green's function for the heterogeneous, absorbing case can be written by analogy as 
\begin{align}
    g_{3D}(\omega,\bx; \bx') &= A(\bx;\bx') \exp{ \Big( i \phi(\bx;\bx') \Big)}.
\end{align}
How the phase $\phi$ and amplitude $A$ are computed in the heterogeneous, absorbing case will be described below in Secs.\ \ref{sec:ray_approx},
\ref{sec:raytracing} and \ref{sec:Discretised-approximation}.

\subsection{Ray-based approximation}
\label{sec:ray_approx}
This section describes, in general terms, how ray theory will be used to compute the terms in the approximate Green's function for heterogeneous absorbing media. 
(See \cite{rullan2018,rullan2020} for a similar approach in the non-absorbing case.)

The numerical implementations are described later, in 
Secs.\ \ref{sec:raytracing} and \ref{sec:Discretised-approximation}.

\subsubsection{Dispersive eikonal equation}
Substituting the Green's function, \eqref{eq:2D_greens} or the 3D equivalent, into the lossy Helmholtz equation, and making the assumption, reasonable for soft tissue, that $\alpha^2 \ll k^2$, leads to the following two equations:
\begin{align} \label{eq:real_parts}
k^2A + \nabla^2A - A\nabla\phi\cdot\nabla\phi = 0,\\
2k\alpha A + 2\nabla A\cdot\nabla\phi + A\nabla^2\phi = 0. \label{eq:imag_parts}
\end{align}
By making the high frequency approximation $| \nabla^2 A / A | \ll k^2$, the first equation leads to the \textit{dispersive eikonal} equation
\begin{align}
\label{eq:dispersive_eikonal}
\nabla \phi \cdot \nabla \phi = k^2.
\end{align}

\subsubsection{Rays}
\label{sec:rays}
In this paper, we solve \eqref{eq:dispersive_eikonal} using the concept of rays. Rays are curves that are perpendicular to surfaces of constant phase, ie.\ they are tangent to the wavevector $\bk$, which satisfies
\begin{align}
\label{eq:phi_grad}
\bk = \nabla \phi,
\end{align}
and therefore \eqref{eq:dispersive_eikonal}. The following definition of a ray will be useful.
\begin{definition} 
A ray connecting an emission point $\bx_e$ to an arbitrary point $\bx \in \Omega $ is defined using $f_{(k,\bx;\bx_e)} (\bx') = 0$. A ray is parameterised by the arc length along it, $s$, by describing it as the line of points $\{ \bx'(s), s\in[0,L_{\text{ray}}]\}$ that satisfy $f_{(k,\bx;\bx_e)}(\bx') = 0$, where $L_{\text{ray}}$ is the physical arc length of the ray, and $s=0$ corresponds to the emission point $\bx_e$. \end{definition}
The numerical procedure for tracing the rays such that they obey \eqref{eq:dispersive_eikonal} will be described in Sec.\ \ref{sec:raytracing}. Once the rays have been computed, the phase and amplitude of the acoustic field can be calculated along the rays.

\subsubsection{Phase along a ray} 
The phase in the homogeneous, non-dispersive, case is \  $\phi_0 (\bx;\bx_e) = k_0|\bx-\bx_e|$, as we saw above. In the heterogeneous, dispersive case, the phase on a ray can be calculated as the line integral 
\begin{align} 
\label{eq:phase}
\phi(\bx;\bx_e) &= \int_{C} \bk \cdot d \bx'  +  \phi(\bx_e),\\
&= \int_{\Omega} k(\bx')  \delta ( f_{(k, \bx; \bx_e)}  )  \ d \bx'  +  \phi(\bx_e). \label{eq:phase2}
\end{align}
where the curve $C$ is the ray joining the start point, $\bx_e$, to the point of interest $\bx$. The second form arises from the fact that the tangent to curve $C$ is everywhere parallel to $\bk$.

\subsubsection{Absorption along a ray} 
The amplitude will be written as two factors,
\begin{align}
\label{eq:A_factors}
A = A_{abs}A_{geom},
\end{align}
where the first accounts for amplitude decay due to absorption losses, and the second for geometric spreading and refractive effects. The absorbing factor is the decay along the ray due to absorption. In the homogeneous but absorbing case we have
\begin{align}
A_{abs}(\bx;\bx_e) = \exp{ \left(- \bk_i\cdot(\bx - \bx_e)  \right) },
\end{align}
which generalises in the case of a non-straight ray to the line integral
\begin{align}
\label{eq:A_abs}
A_{abs}(\bx;\bx_e) &= \exp{ \left(- \int_C \bk_i \cdot d\bx'  \right) }
= \exp{ \left(- \int_C \alpha \left(\frac{\bk}{k}\right) \cdot d \bx'  \right) },\\
&= \exp{\left(  - \int_{\Omega} \alpha(\bx')  \delta( f_{(k, \bx; \bx_e)})  \ d \bx' \right) }.
\end{align}
\subsubsection{Geometric spreading and refractive losses along a ray} \label{sec:geomterical_spreading}
By substituting \eqref{eq:A_factors} into \eqref{eq:imag_parts}, with $A_{abs}$ given by \eqref{eq:A_abs}, gives
\begin{align}
2k\alpha A_{geom} + 2\nabla A_{geom}\cdot\nabla\phi
-2A_{geom}\frac{\alpha }{k}\bk\cdot\nabla\phi
+ A_{geom}\nabla^2\phi = 0.
\end{align}
where we have used 
$\nabla A_{abs} = -(\alpha \bk/k)A_{abs}$.
Using \eqref{eq:phi_grad}, this collapses to
\begin{align}
\label{eq:transport}
\nabla\cdot(A_{geom}^2 \nabla\phi) = 0. 
\end{align}
This \textit{transport equation} describes the effect of geometric spreading (including refractive bending) on the amplitude of the wave. From \eqref{eq:transport}, it is possible to  derive \textit{Green's law}, which relates the amplitude at two different points along the same ray through knowledge of the ray tube area, $S$, at those points \cite{Pierce1981}:
\begin{align}   \label{eq:greens-law}
A_{geom}(\bx) = \Big[ \frac{ (S / \rho c)_{\bx_{ref}}}{(S /\rho c)_{\bx} } \Big]^{1/2} A_{geom} (\bx_{ref}).
\end{align}
How the reference point $\bx_{ref}$ is chosen will be explained below. Sec.\ \ref{sec:Discretised-approximation}, will describe the numerical implementation of the equations in this section.

\subsubsection{Reference point}  \label{reference-point}
The area of ray tube at the emission point vanishes \cite{Cerveny2001}, so, for a finite pressure source $s(t, \bx_e)$, if the emission point $\bx_e$ were chosen as the reference point $\bx_{ref}$ then the pressure amplitude calculated using \eqref{eq:greens-law} would vanish everywhere. To avoid this, the reference point is chosen to be $\bx_\text{ref} = \bx_{e'} \in \mathcal{N}(\bx_e)$ on the ray connecting $\bx_e$ to $\bx$, where $\mathcal{N}(\bx_e)$ is a small neighborhood of the emission point $\bx_e$ and the distance $\lvert \bx_{e'}-\bx_e \rvert$ is a small scalar value. (Here, the reference point is chosen as the first ray point after the emission point.) 

Note that $\Lim{\bx_{e'} \to \bx_e}  S(\bx_{e'}) A_{geom}(\bx_{e'})$ is required to be nonzero and bounded and satisfy \eqref{eq:greens-law} for a given $A_{geom}(x)$ in a homogeneous medium \cite{Cerveny2001}. For this to remain true for the heterogeneous case, the acoustic properties inside $\mathcal{N}(\bx_e)$ are assumed homogeneous.
$A(\bx_{e'})$ is then calculated as the amplitude of the analytic Green's function solution to the wave equation for homogeneous media using \eqref{eq:amplitude-homogenous}.
Correspondingly, given the source function and the emitter position, \eqref{eq:greens-law} can be used to calculate the pressure amplitude at an arbitrary point $\bx$ under an assumption that the area of the ray tubes does not vanish across the medium.

\section{Ray-based inversion accounting for scattering}
\label{sec:optimisation-greens}
In this section, the approach taken for minimising the objective function \eqref{objective-function} in terms of the approximate Green's function introduced in Secs. \ref{sec:greens_functions} and \ref{sec:ray_approx} is explained. The objective function in terms of Green's function in the frequency domain is in the form
\begin{align}  \label{eq:objective-function-greens}
    \mathcal{F}(c) = \frac{1}{2}  \int  \sum_{e,r} \lvert \delta P_{res} \rvert^2 d \omega ,
\end{align}
where $\delta P_{res}$ is the residual, and has components
\begin{align}  \label{res}
\delta p_{res}(c; \omega , r; e )   = p_{(c; \omega , r; e )} - \hat{p}_{( \omega , r; e )},
\end{align}
where $p_{(c; \omega , r; e )} = g(\omega, \bx_r;\bx_e) s(\omega, \bx_e)$. 
\noindent
The minimisation of \eqref{eq:objective-function-greens} involves moving to the minimum of an objective function in steps from some starting point $c^{(0)}$. 
The approach taken in this study is finding the steps $\delta c^{(n)}$ which minimise the Taylor series approximation \cite{Margrave}
\begin{align} \label{eq:objective-function-taylor-series}
    \mathcal{F} ( c^{(n)} + \delta c^{(n)} )  \approx F( c^{(n)} )   + \int \ \frac{ \partial \mathcal{F} }{\partial c^{(n)} (\bx') }  \delta c^{(n)} (\bx') d \bx' +  O\big(\delta (c^{(n)})^2 \big) ,
\end{align}
where $n$ is the iteration number, and the second-order terms and higher are truncated. 
Taking the gradient of \eqref{eq:objective-function-taylor-series} with respect to $c$ yields \cite{Margrave}
\begin{align}  \label{eq:linearised-subproblem}
\begin{split} 
    \frac{ \partial \mathcal{F} \big( c^{(n)} + \delta c^{(n)} \big) }{ \partial c^{(n)} (\bx) } 
    & \approx  \frac{ \partial \mathcal{F} ( c^{(n)} ) }{ \partial c^{(n)} (\bx) } + \int  \    \bigg[   \frac{  \partial^2  \mathcal{F} ( c^{(n)} ) }{\partial \big( c^{(n)} (\bx) \big)  \partial \big( c^{(n)} (\bx') \big)} \delta c^{(n)}(\bx') \bigg] d \bx'   \\
    & = \nabla  \mathcal{F}^{(n)}(\bx) + \int \  H^{(n)} (\bx, \bx') \  \delta c^{(n)}(\bx') d \bx',
    \end{split}
\end{align}

where the first term on the right-hand side is the gradient of $\mathcal{F}$ at iteration $n$, and the second term is the action of the Hessian, $H^{(n)} (\bx, \bx'),$ on a perturbation $\delta c^{(n)}$. Finding the perturbation that makes this right-hand side zero is equivalent to find a search direction using a second-order optimisation approach for minimising $\mathcal{F}$. This leads to a linear minimisation sub-problem, which can be solved by matrix-free approaches, and its solution gives a direct step towards the local minimum of the objective function. Taking this approach for minimising the objective function in \eqref{eq:objective-function-greens} using the usual good approximation of the Hessian described below is called a \textit{Gauss-Newton} algorithm, and requires the solution of a \textit{normal equation} for each sub-problem (see Sec.\ \ref{sec:normal_CG}).

Note that replacing the Hessian with the identity in the above equation will give the steepest descent search direction, which is equivalent to taking a first-order optimisation approach, in which the nonlinear objective function \eqref{eq:objective-function-greens} is minimised using a search direction which uses only the information included in the gradient term in \eqref{eq:linearised-subproblem}. First-order approaches \cite{Beck} are currently in widespread use for full-wave inversion \cite{Wang,Matthews,Matthews1,Liva,Lucka}, because they do not require the additional expense of computing the Hessian. Here, however, the use of a second-order minimisation is the key to incorporating the effects of scattering in the inversion. 
Unlike with some forward models, this does not impose a significant computational burden, because the Hessian matrix can be efficiently computed using the Greens functions already used for computing the gradient. To understand why the scattering effects are included, note that the pressure $p_{(c; \omega, r; e)}$ modelled using the ray approach described in Sec.\ \ref{sec:greens} above only captures the accumulated information along a ray linking emitter $e$ and receiver $r$, ie., the first-arriving signal, and scattering, which allows any point in the medium to contribute to the pressure $p_{(c; \omega, r; e)}$, is neglected. Therefore, using only the information included in the gradient for taking a search direction, scattering effects are neglected. On the other hand, \eqref{eq:linearised-subproblem} can be expressed as a minimisation sub-problem of the form
\begin{align}  \label{eq:linearised-subproblem2}
\delta c^{(n)} \approx \argmin_{\delta c}  \,  \frac{1}{2}  \int  \sum_{e,r}  \lvert   \delta p_{(c^{(n)}; \omega , r; e)} (  \delta c  )   - \delta p_{res}(c^{(n)}; \omega , r; e )    \rvert^2 d \omega,
\end{align}
which seeks to find the search direction $\delta c^{(n)}$ which fits the induced perturbed pressure on the receivers to the residual in a least square sense. The order of dependence of the perturbed pressure $\delta p_{(c^{(n)}; \omega , r; e)}$ on $ \delta c$, which is the same as the order of Taylor series in   \eqref{eq:objective-function-taylor-series}, determines the extent to which scattering is included in \eqref{eq:linearised-subproblem2}.
Here, the dependence of $\delta P$ on $\delta c$ is first order, leading to linear minimisation subproblem, in which $\delta P$ captures single-scattering (which is usually sufficient for soft biological tissue). Therefore, by solving   \eqref{eq:objective-function-greens} through a sequence of linear minimisation problems of the form given in \eqref{eq:linearised-subproblem2}, scattering can be included in the image reconstruction.

\subsection{Fr\'echet derivative}
The Fr\'echet derivative, $ \partial p_{(c^{(n)}; \omega, r; e)} / \partial c^{(n)}( \bx) $ indicates the size of the pressure field perturbation $\delta P$ that will arise from a perturbation in the sound speed $\delta c$. 
The gradient and Hessian can be written in terms of the Fr\'echet derivative, as will be shown below. The form of the Fr\'echet derivative of the forward operator \eqref{eq:forward-operator} will depend on the particular form of the lossy Helmholtz equation used, ie. it will depend on the form of the dispersion relation,
\eqref{eq:general_dispersion_relation}. To retain the generality of the derivation, we first find the Fr\'echet derivative $\partial P / \partial k$ with respect to the wavenumber $k$. The Fr\'echet derivative with respect to the sound speed can be found from that straightforwardly using the chain rule:
\begin{align}
\frac{\partial P}{\partial c} = \frac{\partial P}{\partial k}\frac{\partial k}{\partial c},
\end{align}
where $\partial k/\partial c$ depends on the specific form of the dispersion relation, \eqref{eq:general_dispersion_relation}. 
Recall the lossy Helmholtz equation, \eqref{eq:second-order-lossy-omega}. A perturbation of $k(\bx)$ at a single point $\bx'$ by $\delta k$ will perturb the acoustic field everywhere by $\delta p(\bx)$. Substituting these into \eqref{eq:second-order-lossy-omega} gives
\begin{align} 
\left( (k+\delta k)^2 - \alpha^2 + 2i\alpha(k+\delta k) + \nabla^2  \right) (p + \delta p) = -s,
\end{align}
which, using \eqref{eq:second-order-lossy-omega} and neglecting the products of the perturbations, gives
\begin{align} 
\left( k^2 + \nabla^2  \right) \delta p(\bx) =  - p \, \Upsilon \, \delta k \, \delta(\bx - \bx'),
\end{align}
where $\Upsilon = 2(k + i\alpha)$ and $\Upsilon \, \delta k$ is a complex scattering potential. The solution of this can be written in terms of the Green's function,  \eqref{eq:general-greens-source}, 
so we can write the perturbation of the pressure at any point $\bx$ as
\begin{align} 
\label{eq:pressure_perturbation}
\delta p(\bx;\bx') &= \int g(k,\bx;\bx'')  p(\bx'') \Upsilon(\bx'') \, \delta k \, \delta(\bx'' - \bx') d\bx'', \\
&= g(k,\bx;\bx') \Upsilon(\bx') \delta k(\bx') g(k,\bx'; \bx_e) s(\bx_e),
\end{align}
where the latter step used \eqref{eq:pressure-greens-source}. So the Fr\'echet derivative of the pressure at a point $\bx$ when emitted from a source at $\bx_e$ is 
\begin{align} 
\label{eq:Frechet_general}
\frac{\partial p(\bx;\bx_e)}{\partial k(\bx')} = g(k,\bx;\bx') \Upsilon(\bx') g(k,\bx'; \bx_e) s(\bx_e).
\end{align}
In particular, the Fr\'echet derivative of the pressure at the detector position $\bx_r$ can be computed, for the $n$th iteration, using
\begin{align} 
\frac{\partial p_{(k^{(n)},r;e)}}{\partial k(\bx')} = 
g(k^{(n)},\bx_r;\bx') \Upsilon^{(n)}(\bx') g(k^{(n)},\bx'; \bx_e) s(\bx_e),
\end{align}
This has a clear interpretation, reading from right to left: the source $s$ at point $\bx_e$ is propagated to the point $\bx'$ in the domain
by the Green's function $g(\bx'; \bx_e)$, whereupon it is multiplied by $\Upsilon$ and propagated to the receiver point $\bx_r$ by the Green's function $g(\bx_r;\bx')$. 

\subsubsection{Szabo absorption model}
In order to obtain the Fr\'echet derivative with respect to the sound speed $c$, rather than $k$, we need to use a dispersion relation linking these quantities, ie.\ we need to decide on a particular absorption model.
In Sec.\ \ref{sec:szabo_model}, a form of the dispersion relation for Szabo's absorption model was introduced in \eqref{eq:szabo_dispersion_equation}:
\begin{align}
k = \frac{\omega}{c} + \alpha \tan{(\pi y/2)},
\end{align}
so 
\begin{align} \label{eq:Frechet-derivative}
\begin{split}
\frac{\partial p}{\partial c} &= \frac{\partial p}{\partial k}\frac{\partial k}{\partial c} = -\frac{\omega}{c^2} \frac{\partial p}{\partial k},\\
\frac{\partial p_{(c^{(n)};r,e)}}{\partial c(\bx)} &= g(c^{(n)};\bx_r;\bx)  \Upsilon_c^{(n)}(\bx) g(c^{(n)};\bx; \bx_e) s(\bx_e),
\end{split}
\end{align}
where 
\begin{align}
\Upsilon_c^{(n)}(\bx) = \frac{-2\omega}{c(\bx)^2}\left(\frac{\omega}{c(\bx)} + \alpha(\bx)( \tan{(\pi y/2)} + i ) \right).
\end{align}

\subsection{Gradient and Hessian}
The gradient can be written as the action of the Fr\'echet derivative on the conjugate of the residual $\delta P_{res}^*$ as
\begin{align}  \label{eq:functional-gradient}
    \nabla \mathcal{F}^{(n)}\big( \bx \big) &=  \sum_{e,r}  \int   \text{Re} \left\{ \frac{\partial p_{(c^{(n)}; r; e)}}{ \partial c^{(n)}( \bx) }  \  \delta p_{res}^* (c^{(n)}; r; e)  \right\} d \omega,\\
&=  \sum_{e,r}  \int  \text{Re} \left\{ \Upsilon_{c}^{(n)} (\bx)   \big[  g(c^{(n)}; \bx; \bx_e)  s(\bx_e) \big]  \times  \big[   g(c^{(n)}; \bx_r; \bx)  \  \delta p_{res}^* (c^{(n)}; r, e) \big] \right\} d \omega, \label{eq:gradient-greens2}
\end{align}
where we have used \eqref{eq:Frechet-derivative}.
Similarly, the action of the Fr\'echet derivative on the perturbation $\delta P$ gives a commonly used, and very good, approximation to the action of the Hessian on the sound speed perturbation:
\begin{align}
\label{eq:hessian}
    ( H^{(n)}\delta c^{(n)} ) (\bx) &\approx \sum_{e,r}  \int   \text{Re} \left\{ \frac{\partial p_{(c^{(n)}; r; e)}}{ \partial c^{(n)}( \bx) }  \  \delta p^*_{(c^{(n)}; r; e)}  \right\} d \omega,\\
    &= \sum_{e,r} \int  \  \text{Re} \left\{  \Upsilon_{c}^{(n)} (\bx) \big[ g(c^{(n)};\bx; \bx_e) s(\bx_e) \big]  \times \big[  g(c^{(n)}; \bx_r; \bx)  \  \delta p^*_{(c^{(n)}; r; e)}  \big] \right\} d\omega,\label{eq:hessian2}
\end{align}
where again we have used \eqref{eq:Frechet-derivative}. In both of these integrals, all the factors of the integrand depend on frequency $\omega$. The square brackets here and above just indicate a way of interpreting this formulation of the gradient as the product of a forward propagation of the source by the Green's function, and backwards propagating of the residual by its adjoint. In other words, in \eqref{eq:Frechet-derivative}, $ g(\bx_r; \bx) $ acts as a forward propagator from $\bx$ to $\bx_r$, whereas in Eqs.\ \eqref{eq:gradient-greens2} and \eqref{eq:hessian2}, it acts as a backward propagator from $\bx_r$ to $\bx$. Note that the Green's function obeys the reciprocity relation $g(\bx;\bx') = g(\bx';\bx)$.

We emphasise here again the reason for using the Hessian in the inversion. The residual, $\delta p_{res}$, which appears in the expression for the gradient, \eqref{eq:gradient-greens2}, depends only on the points along a ray linking $e$ and $r$, and therefore it includes only refractive effects. However, the pressure perturbation $\delta p$ which appears in the expression for the Hessian, \eqref{eq:hessian2}, can be affected by any point in the medium, because it is obtained from the action of the Fre\'chet derivative, on the sound speed perturbation, $\delta c$, which can vary anywhere. Therefore, solving the minimisation subproblem \eqref{eq:linearised-subproblem2}, which gives the normal equation \eqref{eq:linearised-subproblem}, can account for single scattering in a way that the ray-based gradient cannot. This is the reason why, when using a ray-based forward model, a Gauss-Newton approach to the inversion step, rather than a purely gradient-based approach, is critical for obtaining good reconstructions.

\subsection{Normal equation and update computation}
\label{sec:normal_CG}
Having defined the gradient and Hessian, this section explains how they are included in \eqref{eq:linearised-subproblem} and a Gauss-Newton step direction is calculated. In this case, the linear subproblem, \eqref{eq:linearised-subproblem}, takes the form of the normal equation, which is solved using a \textit{conjugate gradient} (CG) algorithm. The procedure is outlined below.
\begin{algorithm}
	\caption {\textit{Conjugate gradient} (CG) algorithm for computing update}
	\label{alg:conjugate-gradient}
	\begin{algorithmic}[1] \label{hes}
		\State \textbf{input:} $ c^{(n)}, \ \mathcal{F}( c^{(n)} ), l_\text{max} $
		\State \textbf{initialise:} $ l=0, \ \delta c_l^{(n)}=0$
		\State $ r_l \leftarrow \nabla \mathcal{F}^{(n)} $ using \eqref{eq:gradient-greens2} 
		 \Comment{Initialise the cg (inner) residual}
		\State $ d_l =  r_l  $
		\Comment{Initialise the cg (inner) step direction}
		\While {$l < l_{\text{max}} $}
		    \Comment{Iterate for a fixed number of iterations}
		    \State {$z_l \leftarrow H^{(n)} d_l$}
		    \Comment{Update the action of the Hessian on the cg step direction using \eqref{eq:hessian}}
            \State   $\alpha_l  \leftarrow \frac{r_l^T r_l}{d_l^T z_l }$ 
            \Comment{Update the cg step size 1}
		    \State   $\delta c^{(n)}_{l+1} \leftarrow \delta c_l^{(n)} + \alpha_l d_l $
			\Comment{Update the sound speed perturbation}
		    \State  $ r_{l+1} \leftarrow r_l - \alpha_l z_l $ 
		     \Comment{Update the cg residual}
		    \State   $\beta_l  \leftarrow \frac{r_{l+1}^T r_{l+1} } {r_l^T r_l }$ 
            \Comment{Update the cg step size 2}
            \State $d_{l+1} \leftarrow r_{l+1} + \beta_l d_l $
            \Comment{Update the step direction}
            \State $ l   \leftarrow l+1 $
            \Comment{Increment the counter}
		    \EndWhile
		\State \textbf{output:}  $ c^{(n+1)} \leftarrow  c^{(n)} + \delta c_*^{(n)} $
		    \Comment{Return the optimal sound speed}
	\end{algorithmic}
\end{algorithm}

In Algorithm \ref{alg:conjugate-gradient}, $l$ is the number of CG iterations, and the CG algorithm is terminated after a fixed $l_\text{max}$ number of iterations. Also, the Green's functions included in $\nabla F^{(n)}$ and $H^{(n)}$ are dependent on the chosen angular frequencies $\omega$. Algorithm \ref{alg:conjugate-gradient} is iteratively implemented by choosing a range including small frequencies for $c^{(0)}$, and increasing the frequency range for the next linearised subproblems.

\section{Numerical Ray Tracing}
\label{sec:raytracing}
This section, and the following section, describe the numerical implementation of the method described above for implementing UST. It is therefore worth taking stock at this point and recapping the essential principles. We aim to reconstruct an image of the sound speed by minimising the objective function $\mathcal{F}$ in \eqref{eq:objective-function-greens} through an iterative formation of the linear subproblem \eqref{eq:linearised-subproblem}, and solving it by means of a CG algorithm given in Algorithm \ref{alg:conjugate-gradient}. The first term in \eqref{eq:linearised-subproblem} is the functional gradient $\nabla \mathcal{F}$, and is calculated using \eqref{eq:gradient-greens2}. The second term is an action of the Hessian matrix on a perturbation to the sound speed, and is calculated using \eqref{eq:hessian2} and \eqref{eq:pressure_perturbation}. These formulae are functions of the Green's functions $g(c^{(n)}; \omega, \bx; \bx_e) $ and $g(c^{(n)}; \omega, \bx_r; \bx) $. The phase in the Green's function is computed from Eqs.\ \eqref{eq:phase}, and the amplitude $A$ is computed from Eqs. \eqref{eq:A_abs} and \eqref{eq:greens-law}. All these formulae rely on the rays introduced in Sec.\ \ref{sec:rays}. How those rays are traced in practice is described in this section.

\subsection{Ray equations}
In Sec.\ \ref{sec:ray_approx}, as a means of solving the dispersive eikonal equation \eqref{eq:dispersive_eikonal}, rays with direction vector $\bk(\bx) = \nabla \phi(\bx)$ were introduced. 
It can be straightforwardly shown that the unit vector $d\bx/ds$, describing how the position vector $\bx(s)$ changes with distance along the ray $s$, and the direction vector $\bk$ must satisfy the coupled \textit{ray equations}
\begin{align}
\label{eq:ray_equations}
\frac{d\bx}{ds} = \frac{\bk}{k}, \quad \frac{d\bk}{ds} = \nabla k.
\end{align}
In UST, the rays must always start at an emission point $\bx_e$ and end at a receiver point $\bx_r$. Computing rays to link two points like this is known as \textit{two-point} ray tracing \cite{Anderson2}. In this context it is instructive to consider Fermat's principle \cite{Pierce1981, Holm}, which states that the path between two points taken by a ray makes the acoustic length stationary under variations in a family of nearby paths. 
\begin{definition} 
\label{def:acoustic_length}
The trajectory of a ray passing from point $p_1$ to $p_2$ is the minimum for the variational integral
\begin{align}
\label{stationary_action}
L_{(p_2,p_1)} = \int_{p_1}^{p_2} k\big( \bx(s) \big) \ ds,
\end{align}
where $ds = |d\bx(s)|$ is an infinitesimal distance along the ray, and $L_{(p_2,p_1)}$ is the acoustic length between points $p_1$ and $p_2$.
\end{definition}
This can be interpreted in a Lagrangian sense as minimising the integral of a Lagrangian $\mathcal{L} = T - V$ in order to find the path taken by a particle with kinetic and potential energies given by
\begin{align}
T = k(\bx) \dot{\bx}\cdot\dot{\bx} /2,\quad V = -k(\bx)/2,
\end{align}
where $\dot{\bx}\equiv d\bx/ds$. In other words, the ray is the path of a particle with velocity $\dot{\bx}$ and (spatially-varying) mass $k$ in a potential field given by $-k/2$. Equivalently, the ray equations could be derived from the Hamiltonian $H(\bx,\bk) = T + V$. Interestingly, many different Hamiltonian systems can be defined whose bicharacteristic curves can be associated with the rays as defined here \cite{Cerveny2001, Chapman, Rawlinson2,Engquist,Runborg, rullan2018}. In practice, Eqs.\ \eqref{eq:ray_equations} are solved using a numerical ray-tracing algorithm combined with \textit{ray-linking} techniques, as described below.

\subsection{Ray tracing with Heun's method}  
Several algorithms for computing rays that satisfy Eqs.\ \eqref{eq:ray_equations} are described in the literature, see eg. \cite{Johnson, Anderson2, Bold, Anderson3, Cerveny2001, Runborg, Rawlinson2}. Specifically, this system of coupled equations has been solved using a second-order variant of \textit{Runge-Kutta} method \cite{Butcher}. This approach is taken in this paper.

Consider a ray, and for a point $m$ along the ray let us define $ y_m = \left[ \bx_m, \bkappa_m \right]^T \in \mathbb{R}^{2d}$. Using a second-order variant of Runge-Kutta method, also known as \textit{Heun's} method, $y_{m+1}$ can be approximated in the form
\begin{align}  \label{Runge-Kutta}
    y_{m+1} = y_m + \frac{ \Delta s}{2} (q(y_m) + q(\hat{y}_{m+1})).
\end{align}
where
\begin{align}  \label{update_variables}
q(y) =
\begin{bmatrix}
 q_{\bx}\\
 q_{\bkappa}
\end{bmatrix}
=
\begin{bmatrix}     
 \bkappa / k(\bx)\\
 \nabla k (\bx) 
\end{bmatrix},
\end{align}
and
\begin{align}   
    \hat{y}_{m+1} = y_m + \Delta s \, q(y_m).
\end{align}
In addition, the step directions are normalised such that $|\bkappa| = k$, which ensures $\dot{\bx}\cdot\dot{\bx} = 1$. The ray tracing algorithm using Heun's method is outlined in Algorithm \ref{alg:ray_tracing_mixed_step}. 

\begin{algorithm}
    \caption{Ray tracing using Heun's method}
    \label{alg:ray_tracing_mixed_step}
    \begin{algorithmic}[1]
        \State \textbf{input:} $\bx_{e}$, $k:= k(\bx)$
            \Comment{Input initial ray position and wavenumber}
            \State \textbf{initialise:} $\bx = \bx_{\bar{e}}$, $\bkappa = \bkappa_{\bar{e}}$           \Comment{Set initial ray position and direction}
        \While { $\bx(s) \ \text{is inside} \ \Omega$ } 
        \State $\bkappa \leftarrow k \   \bkappa / \mid \bkappa \mid $
           \Comment{Normalise the ray direction}
            \State $q_{\bx} = \bkappa / k$
            \Comment{Compute the update variables using \eqref{update_variables}   }  
            \State $q_{\bkappa} = \nabla k(\bx)$  
           \State $\bkappa' \leftarrow \bkappa + \Delta s  \ q_{\bkappa} $
      \Comment{Update the auxiliary ray direction}
        \State $k' \leftarrow k (\bx + \Delta s \  q_{\bx})$
        \Comment{Update the auxiliary wavenumber}
               \State $\bkappa' \leftarrow   k'  \bkappa' / \mid \bkappa' \mid $
        \Comment{Normalise the auxiliary ray direction}  
             \State $q_{\bx}' = \bkappa' / k' $
              \Comment{Compute the auxiliary update variables using \eqref{update_variables}  }
             \State $q_{\bkappa}' = \nabla k(\bx + \Delta s  \ q_{\bx})$   
             \State $\bx \leftarrow \bx + (\Delta s/2) (q_{\bx} + q_{\bx}')$
                \Comment{Update the ray position using \eqref{Runge-Kutta}}  
            \State $\bkappa \leftarrow \bkappa + (\Delta s/2) (q_{\bkappa} + q_{\bkappa}')$
                \Comment{Update the ray direction using \eqref{Runge-Kutta}}  
        \EndWhile
    \end{algorithmic}
\end{algorithm}

\newpage

\subsection{Grid-to-ray interpolation}
\label{subsec:Grid-to-ray_interpolation}
The wavenumber field $k(\bx)$ was represented on a discretised mesh of $N_n$ points $x_i$. A rectilinear grid was used with grid points indexed with the multi-index $i = \big( i^1 , i^2  \big) \in \left\{ 1, ... , N_n^1 \right\}  \times \left\{ 1, ... , N_n^2 \right\}  $ with $N_n = \prod_{j=1}^2 N_n^j $ and an equal grid spacing $\Delta x$ along all Cartesian coordinates $j$. (Recall that our study is restricted to $d=2$.) Also, $x_{i^j}$ is used to indicate the position of grid point $i$ along Cartesian coordinate $j$. The points along a ray can lie on any arbitrary points in $\Omega$, and are not restricted to the grid points. Therefore, an interpolation from the grid to the rays must be performed to find the approximate values for $k$ and  $\nabla k$ used in the ray tracing algorithms \cite{Anderson2,Denis}. Here, the wavenumber function $k(\bx)$ is represented with a set of B-spline functions, which therefore gives continuous values for the directional gradients. Following \cite{Virieux,Denis}, the control points for the B-spline function are chosen on the grid points $Q(x_{i^1},x_{i^2})$, and the B-spline function is defined using
\begin{align}  \label{eq:bspline-map}
    k (\bx) \approx \hat{k}(\bx) =  \sum_{q^1=0}^3 \sum_{q^2=0}^3 Y_{q^1}(u) Y_{q^2}(v)  Q(x_{i^1+q^1-1},x_{i^2+q^2-1}).
\end{align}
Here, $u$ and $v$ are defined as 
\begin{align}
\begin{split}
u(\bx) = 
    \begin{cases}
      \frac{x^1- x_{i^1}}{\Delta x}, & \text{if}\ x_{{(i-1)}^1} < x^1 <  x_{{(i+1)}^1} \\
      0, & \text{otherwise}
    \end{cases}
 \\
v(\bx) = 
    \begin{cases}
      \frac{x^2- x_{i^2}}{\Delta x}, & \text{if}\ x_{{(i-1)}^2} < x^2 <  x_{{(i+1)}^2} \\
      0, & \text{otherwise},
    \end{cases}
\end{split}
\end{align}
where $x_{{(i\pm1)}^j}$ denotes the two grid points adjacent to the grid point $x_i=(x_{i^1},x_{i^2})$ along the Cartesian coordinate $j$. Here, the polynomials $Y_m(u)$ and $Y_n(v)$ satisfy 
\begin{align} \label{eq:bspline-polynomials}
\begin{split} 
\begin{bmatrix}
Y_0(u)\\
Y_1(u)\\ 
Y_2(u)\\
Y_3(u)
\end{bmatrix}=
\frac{1}{6}
\begin{bmatrix}
-1 &3 &-3 &1\\
3 &-6 &0 &4 \\ 
-3 &3 &3 &1\\
1 &0 &0 &0
\end{bmatrix}
\begin{bmatrix}
u^3\\
u^2\\ 
u\\
1
\end{bmatrix}.
\end{split}
\end{align}
Using Eqs. \eqref{eq:bspline-map} and \eqref{eq:bspline-polynomials}, the approximated wavenumber field is $\mathcal{C}^2$ continuous \cite{Virieux,Denis}. Therefore, the components of the directional gradients can be analytically calculated from \eqref{eq:bspline-polynomials}.

\subsection{Ray-linking}  \label{sec:raylinking}
The rays initialised from an emitter at point $\bx_e$ are used for two tasks, (1) computing the Green's function at the reception points $\bx_r$, i.e., $g(c^{(n)}; \omega , \bx_r; \bx_e)$, which is used for calculating the residual $\delta P_{res}$, and (2) computing the Green's function at arbitrary points within the medium $\bx$ from excitation $e$, ie., $g(c^{(n)}; \omega , \bx; \bx_e) $.
Tracing individual rays from $\bx_e$ to every point in the domain separately would be much too expensive and many of the rays would be redundant. There are therefore two choices: (a) send rays from $\bx_e$ at evenly-spaced initial angles out across the whole domain with no specific end-point in mind, and interpolate from these rays onto $\bx$ and $\bx_r$ as required, or (b) trace rays from $\bx_e$ to $\bx_r$ using ray-linking and then, using those rays, interpolate onto the points $\bx$ as required. The downsides of option (a) are that there will be interpolation error in the estimate of the pressure at the reception point $\bx_r$. Therefore, the spacing for the initial angles should be sufficiently small. More importantly, the acoustic length between two points on such a ray is less likely to be a global minimum \cite{Anderson2}. 
With option (b), on the other hand, ray-linking seeks to find a ray trajectory which provides the stationary point for the acoustic length between $\bx_e$ and $\bx_r$ by enforcing a boundary condition on the rays' path such that the ray initialised from an emitter $e$ is intercepted by the receiver $r$ after travelling across the medium. 
The main downside of (b) is that \textit{per ray} it is more expensive than (a) because the ray-linking is iterative. However, there are several advantages to using ray-linking. First, 
because of the imposition of the boundary condition on the ray, the approach is more stable in the following sense: the trajectory of the ray for option (a) doesn't link $e$ to $r$ but to a point different from $r$, and the phase and amplitude difference between this point and $r$ may be too large to correct using a simple interpolation, resulting in errors in the phase and amplitude.
Second, the linked rays can be re-used when calculating the adjoint field, as they link directly to the receiver points from which the adjoint sources emanate. Finally, the number of rays is fixed at $N_rN_e$, whereas in (a) it could be greater, depending on the criterion for choosing the initial angular spacing. 

The most popular ray-linking approaches for 2D media are based on solving the associated one-dimensional inverse problem using the \textit{Regula falsi} or \textit{Secant} approaches \cite {Anderson2,Cerveny2001}. The \textit{Regula falsi} approach will be used for validation of the ray tracing using a digital breast phantom, and a \textit{Secant} approach will be used iteratively for solving the optmisation problem in a way in which the linked rays used for solving the linearised problem $n$ will be used as an initial guess for the ray linking for the next subproblem $n+1$. (See \cite{Javaherian} for an extension of ray-linking to 3D for a UST inversion approach using time-of-flights.) 

\subsection{Ray coordinates} 
In the general case, the coordinates of the ray are given by two parameters: one specifying the initial direction of the ray and another monotonic parameter along the ray \cite{Cerveny2001}. These parameters can be chosen in different ways. For our 2D case, the ray coordinates will be defined as the initial direction of the ray in polar coordinates, which is also called \textit{radiation angle}, and the arc length $s$ along the ray, with $s=s_0$ matching the emission point, and monotonically growing along the ray.

\begin{definition}  \label{def:ray-parameterisation}
The trajectory of a ray linking an emission point $e$ to a reception point $r$ is defined by the points with arc lengths $s_m, \ m \in \left\{ 0,...,M_{(e,r)} \right\} $, where the number of the sampled points along the ray is $M_{(e,r)}+1$. Therefore, the points are initialised from $s_0$ with $\bx(s_0) := \bx_e$, and are terminated at the point $s_{M_{(e,r)}}$ with $\bx(s_{M_{(e,r)}}) :=\bx_r$, the position of receiver $r$. The points $s_m$ satisfy
\begin{align}  
\label{eq:ray_points_sm}
\begin{split}
s_m = 
\begin{cases}
m \Delta s,                                & m\in\left\{0,..., M_{(e,r)} -1 \right\}\\
\left( m-1 \right)  \Delta s + \Delta s',   & m= M_{(e,r)}.
\end{cases}
\end{split}
\end{align}
Here, the second line is used in order to indicate that the last point of the ray must be matched to the reception point $r$, and thus $\Delta s' = s_{M_{(e,r)}} - s_{M_{(e,r)-1}}$ with $ \Delta s'  \leqslant  \Delta s$ \cite{Javaherian}. 
\end{definition}

\begin{definition}  \label{def:raylinked-coordinates}
For each excitation $e$, the pressure field is approximated on a set of linked rays  $f_{(k,\bx_r;\bx_e)} = 0, \ e \in \left\{ 1,...,N_e \right\},  \  r \in \left\{ 1,...,N_r \right\}$. These rays are parameterised in space using $\bx(s_m,\theta_{(r,e)})$, which denotes the position on the arc length $s$ of the point $m$ along a ray linking the emission position $\bx_e$ to the reception position $\bx_r$. Also, the polar initial direction of this ray is indicated by $\theta_{(r,e)}$. 
\end{definition}

Having defined the rays linking all emitter-receiver pairs, the sampled points along these rays are used for a spatial parameterisation of the medium and for a numerical implementation of the method approximating the acoustic field as described in Sec. \ref{sec:greens}. The discretised version of this theory is the topic of the next section.



\section{Discretisation of approximate Green's function} \label{sec:Discretised-approximation}

This section describes how the model outlined in Sec.\ \ref{sec:greens} is discretised for implementation. The Green's function, $g(\omega, \bx; \bx_e) $, is discretised at the sampled points along the rays linking the emitter $e$ to all receivers $r$ using the coordinates defined in Definition \ref{def:raylinked-coordinates}
\begin{align}  \label{eq:greens-heterogeneous-sampled}
\begin{split}
    & g \big( \omega; \bx( s_m, \theta_{(r,e)} ) \big) \approx  \\
    & A \big( \bx(  s_m;  \theta_{(r,e)}  ) \exp{ \Big( \iu \big(  \phi   \big( \bx( s_m, \theta_{(r,e)})  \big) + \pi/ 4  \big)   \Big)  },
    \end{split}
\end{align}
where $ \phi  \big( \bx( s_m, \theta_{(r,e)})  \big) $ and $ A  \big( \bx( s_m, \theta_{(r,e)})  \big)$ are the phase and amplitude on point $m$ on the ray linking emitter $e$ to receiver $r$, respectively.

\subsection{Acoustic absorption}
The acoustic absorption on the ray is computed 
in the form
\begin{align}  \label{eq:amplitude-abs-discretised}
    A_{abs}\big( \bx( s_m, \theta_{(r,e)})  \big) =   \exp{ \left( -\int_{s_0}^{s_m} \alpha \big( \bx(s_m, \theta_{(r,e)})\big) \ ds   \right) },
\end{align}
where $A_{abs} \big( \bx(s_0, \theta_{(r,e)})\big) = 1 $.

\subsection{Geometrical spreading}
\label{sec:Discretised_geometerical_spreading}

The effect of geometric spreading on the amplitude of Green's function is described in \eqref{eq:greens-law}, which relates the amplitudes at two points along a ray by considering how the area of the \textit{ray tube} around the ray changes between the two points, as described in Sec.\ \ref{sec:geomterical_spreading}. To compute this numerically requires us to estimate the rate at which two closely-spaced rays diverge \cite{rullan2018}, for which we turn to the concept of the \textit{ray Jacobian} \cite{Cerveny2001, rullan2018}.

First, note that the relations between the ray coordinates, $\gamma = \left[ \gamma_1, \gamma_2 \right]^T = \left[ \theta_{(r,e)}, s_m \right]^T$ and the general Cartesian coordinates, $\bx = \left[ x_1, x_2 \right]^T$, follows 
\begin{align}
  Q_{ii'}(s_m, \theta_{(r,e)})= \frac{\partial x_{i}}{\partial \gamma_{i'}}(s_m, \theta_{(r,e)}),
\end{align}
where $Q$ is called the \textit{transformation matrix}, and satisfies 
\begin{align}
    d \bx(s_m, \theta_{(r,e)}) =  Q (s_m, \theta_{(r,e)}) d\gamma.
\end{align}
The Jacobian $J$ on the point $m$ along the ray initialised by the angle $\theta_{(r,e)}$ then satisfies
\begin{align}
    \label{eq:Qdet}
    J(s_m, \theta_{(r,e)}) = \text{det} \  Q (s_m, \theta_{(r,e)}),
\end{align}
where $\text{det}$ denotes the determinant. This ray Jacobian is closely connected to the density of the ray field, which is described by the cross section area of ray tube. Here, the ray tube for the ray $\theta_{(r,e)}$ is defined as a family of rays with an initial angle in the interval $\big(\theta_{(r,e)}-\Delta \theta , \ \theta_{(r,e)}+\Delta \theta  \big)$.
We can estimate the Jacobian $J$ as the determinant of a matrix which includes the derivative of a ray's trajectory with respect to the ray coordinates using finite differences \cite{rullan2018}
\begin{align}
    \label{eq:Jacobian}
    J \big( s_m, \theta_{(r,e)}\big)  = \text{det} \bigg( \Big[  \frac{\partial \bx }{\partial \theta}\big(s_m, \theta_{(r,e)} \big), \ \frac{\partial \bx}{\partial s} \big(s_m, \theta_{(r,e)} \big)  \Big]^T  \bigg).
\end{align}
The derivative of the ray's trajectory with respect to the arc length $s$ can be approximated using finite differences in the form
\begin{align}  \label{eq:derivative_s}
\frac{\partial \bx}{\partial s}  \big(s_m, \theta_{(r,e)} \big) \approx \frac{\bx\big(s_{m+1} , \theta_{(r,e)} \big) -\bx  \big( s_{m-1}, \theta_{(r,e)} \big) }{2\Delta s},
\end{align}
where $\Delta s $ is the user-defined step size used for calculation of the trajectory of the ray (cf. Sec. \ref{sec:raytracing}). In the same way, the derivative with respect to the polar initial direction can be numerically approximated using 
\begin{align}  \label{eq:derivative_angle-reference}
\frac{\partial \bx}{\partial \theta}
\big(s_m, \theta_{(r,e)} \big) \approx \frac{\bx\big(s_m, \theta_{(r,e)} + \Delta \theta \big) -\bx\big(s_m, \theta_{(r,e) } - \Delta \theta \big) }{2\Delta \theta },
\end{align} 
where $\Delta \theta$ is a user-defined perturbation in the initial angle of the linked ray, and is here fixed. The above equation requires tracing two additional auxiliary rays with angles $\theta_{(r,e)} \pm \Delta \theta$ for each emission-reception pair after ray linking \cite{Cerveny2001}.
When the reception points are sufficiently close together the computation of auxiliary rays can be avoided by modifying \eqref{eq:derivative_angle-reference} into the form
\begin{align}  \label{eq:derivative_angle-raylinked}
\frac{\partial \bx}{\partial \theta}
\big(s_m, \theta_{(r,e)} \big) \approx \frac{\bx\big(s_m, \theta_{(r+1,e)}   \big) -\bx\big(s_m, \theta_{(r-1,e) } \big) }{  \theta_{(r+1,e)}-\theta_{(r-1,e)} }.
\end{align} 
Here, $\bx \big(s_m, \theta_{(r\pm 1, e)} \big)$ denotes the position of the two rays linking the emission position $\bx_e$ to the two nearest reception points $\bx_{r \pm 1}$. (Note that a rotational indexing must be used for the reception points for a circular detection geometry, ie.\ the receivers adjacent to $r=N_r$ are $r=\left\{N_r-1, 1\right\}$.)
Using \eqref{eq:derivative_angle-raylinked}, the two auxiliary rays used in \eqref{eq:derivative_angle-reference} are replaced by two adjacent linked rays, which have already been calculated. 
Considering \eqref{eq:greens-law}, the geometrical attenuation satisfies
\begin{align}  \label{eq:amplitude-geom-discretised}
A_{ \text{geom} } \big(  \bx( s_m ; \theta_{(r,e)} )  \big) = \bigg[ \frac{ c \big(\bx(s_m,\theta_{(r,e)})\big)  }{c \big(\bx(s_1,\theta_{(r,e)})\big)} \frac{J(s_1,\theta_{(r,e)})}
 {J (s_m,\theta_{(r,e)} )} \bigg]^{1/2},
\end{align}
where  $J\big(s_1,\theta_{(r,e)}\big)$ is the Jacobian of the ray on the reference point, which is chosen as the second point along the ray, using an assumption that a neighborhood of the emission point with a radius greater than ray spacing $\Delta s$ is acoustically homogeneous. It is reminded that the acoustic density has been assumed homogeneous in this study.

\subsection{Phase}
The accumulated phase is discretised in the form
\begin{align}  \label{eq:phase-discretised}
    \phi   \big( \bx( s_m, \theta_{(r,e)})  \big) = \int_{s_0}^{s_m} k\big( \bx(s_m, \theta_{(r,e)})\big) \ ds   + \phi  \big( \bx( s_0, \theta_{(r,e)})  \big) + \frac{\pi}{2}  K(s_m, \theta{(r,e)} ) .
\end{align}
where $\bx( s_0, \theta_{(r,e)}) := \bx_e $ is the emission point. Also, $K(s_m, \theta{(r,e)} )$ is the cumulative times the sign of the ray Jacobian along the ray has been changed. Points on which the ray Jacobian changes sign are called \textit{caustics}, and will lead to a $\pi/2$ shift in the phase \cite{Cerveny2001}.

\section{Numerical Results}  
\label{sec:num-results}
This section will describe numerical experiments demonstrating the effectiveness of the proposed ray-based inversion approach for a quantitative reconstruction of the sound speed distribution of an object from the measured pressure time series. In Sec. \ \ref{sec:data-simulation}, the details for simulation of ultrasound pressure time series are explained. In Sec. \ref{sec:results-greens}, the Green's function solution to the Szabo's second-order wave equation \eqref{eq:second-order-lossy} is compared to a full-wave solution using a \textit{k-space pseudo-spectral} method \cite{Cox,Treeby2}. This numerical wave solver is freely available in the open-source \textit{k-Wave} toolbox \cite{kWave,Treeby}. The reconstructed images will be presented in Sec. \ref{image_reconstruction}.

\subsection{Data simulation} \label{sec:data-simulation}
\subsubsection{Imaging system and breast phantom}
\begin{figure} 
    \centering
	\subfigure[]{\includegraphics[width=0.45\textwidth]{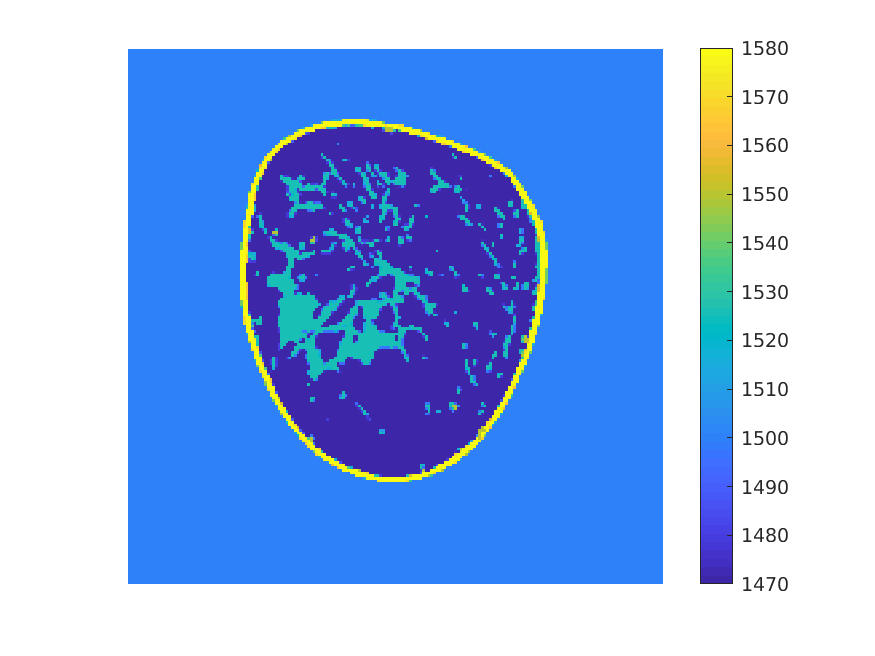}
     \label{fig:phantom_sound_speed}  }
      \subfigure[]{\includegraphics[width=0.438\textwidth]{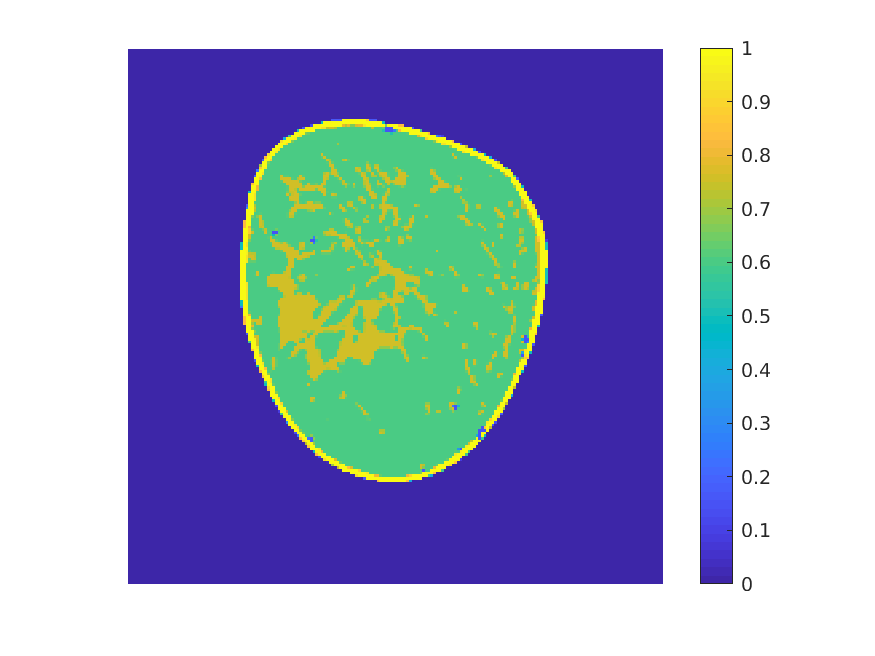}
     \label{fig:phantom_absorption_coefficient}   }
	\caption{Phantom: (a) sound speed [m/s] (b) 
	absorption coefficient [$\text{dB} \text{MHz}^{-y}\text{cm}^{-1}$],
    and a homogeneous power law exponent $y=1.4$.}
\end{figure}

Ultrasound tomography data was simulated for an imaging system which consists of 64 emitters and 256 receivers uniformly distributed along a circle with radius $R=9.5$cm. A horizontal slice of a 3D digital phantom mimicking the acoustic properties of the breast was used in this study. This digital phantom is freely available \cite{Lou}. The sound speed was set to a range $1470-1580$m/s, and the absorption coefficient $\alpha_0$ was set to a range $0-1 \ \text{dB} \text{MHz}^{-y} \text{cm}^{-1} $, and the power law exponent $y$ was set $1.4$. The sound speed and absorption phantoms are shown in Figs. \ref{fig:phantom_sound_speed} and \ref{fig:phantom_absorption_coefficient}, respectively.
The computational grid consisted of $670 \times 670$ grid points with position $\big[-10.05, +10.02\big] \times \big[-10.05, +10.02\big] \text{cm}^2$ and a grid spacing of $ \Delta x = 3 \times 10^{-2}$cm along all the Cartesian coordinates. With this sound speed distribution and grid spacing, the maximum frequency  supported by the grid, $f_{max}$, was $2.45$MHz.

\subsubsection{Simulating time series data}
\begin{figure} 
   \centering
	\subfigure[]{\includegraphics[width=0.45\textwidth]{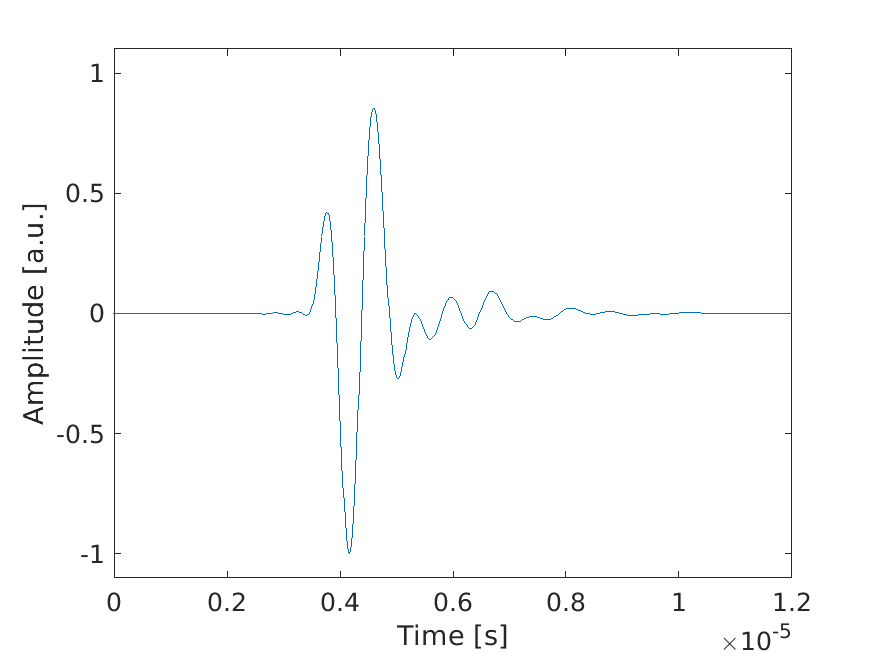}
   \label{fig:pulse_time}   }
   \subfigure[]{\includegraphics[width=0.45\textwidth]{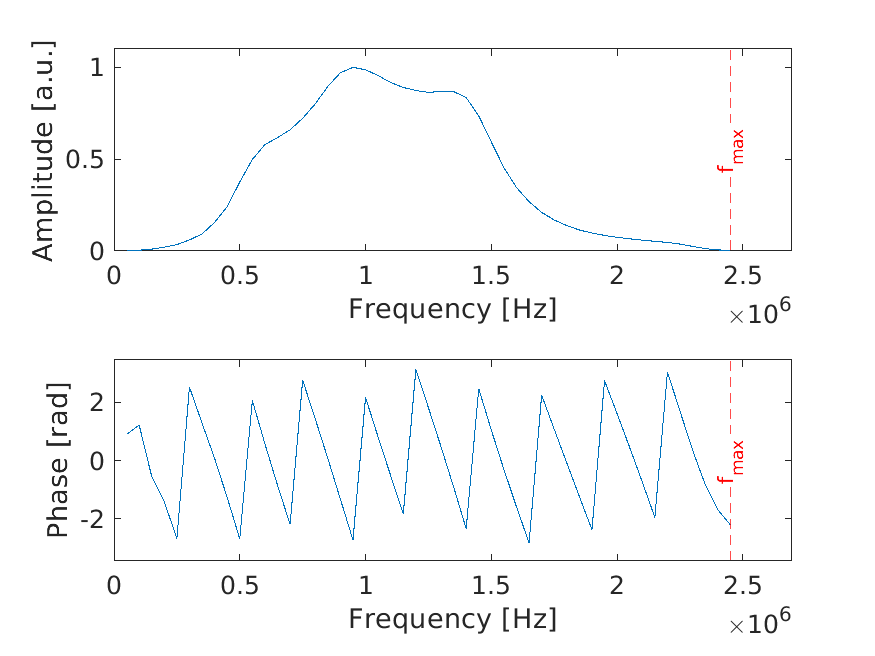}
    \label{fig:pulse_frequency}  }
	\caption{Excitation pulse used for all emitters: (a) time domain, (b) frequency domain: normalised amplitude and phase. $f_{max}$ indicates the maximum frequency supported by the grid used in the k-wave simulations.}\label{fig:excitation-pulse}
\end{figure}

A \textit{k-space pseudo-spectral} method (k-Wave) was used for simulation of the acoustic pressure time series data \cite{kWave,Treeby}. 
The emitters and receivers were assumed as points (not necessarily lying on grid points), and the interpolation of the pressure field from the grid to these transducers and vice versa was performed using the \textit{off-grid} toolbox \cite{Wise}. To simulate the data, each emitter was individually driven by an excitation pulse, and the acoustic pressure time series induced at the receivers were recorded simultaneously. This was repeated for each emitter. The pressure time series were recorded at 8621 time points with a sampling rate of $52.8$ MHz ($1.8939 \times 10^{-8} \ s$ time spacing). Additive white Gaussian noise was added to the simulated pressure time series to give a $40$ dB signal-to-noise ratio of the peak amplitudes. 
The output of physical ultrasound transducers covers typically a finite and quite limited bandwidth, and the generated pressure field tends to be more directional at high frequencies. In designing transducers for an imaging system, therefore, a trade-off must be made between the range of frequencies in the excitation pulse and the directionality of the detectors. Although here the transducers are assumed omnidirectional at all frequencies, the frequency bandwidth of the simulated excitation pulse was nevertheless limited. Fig.\ \ref{fig:pulse_time} shows the normalised amplitude of the excitation pulse (pressure source) in the time domain, and Fig. \ref{fig:pulse_frequency} shows the normalised amplitude and phase components of the excitation pulse in frequency domain, respectively. This signal is used as the pressure source for all excitations. Note that for practical experiments, the induced pressure source may not be the same for different emitters, because it depends on the properties of the emission elements, and a calibration step may be necessary.

\begin{figure}
   \centering
	\subfigure[]{\includegraphics[width=0.45\textwidth]{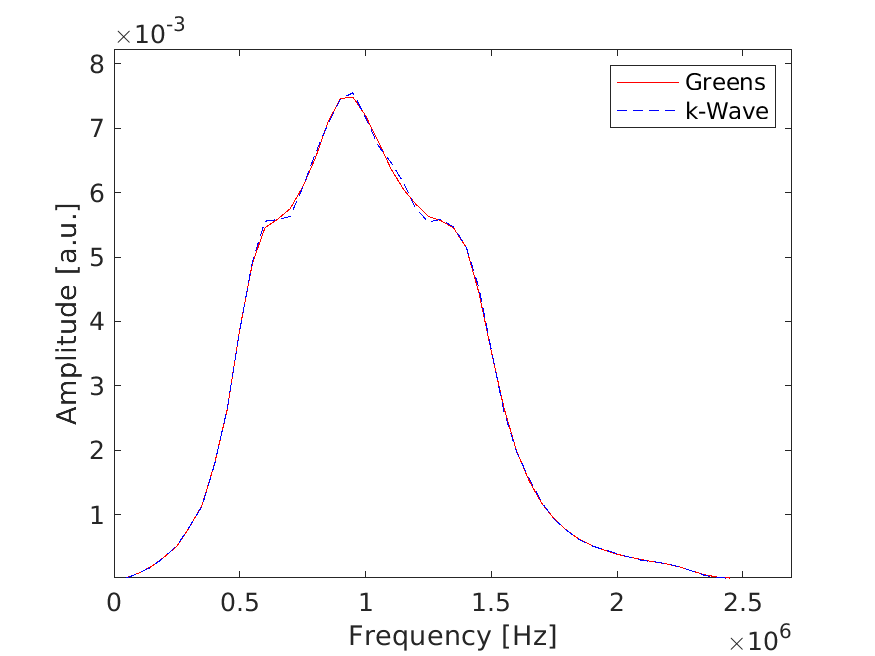}
   \label{fig:signal_water_amplitude} }
   \subfigure[]{\includegraphics[width=0.45\textwidth]{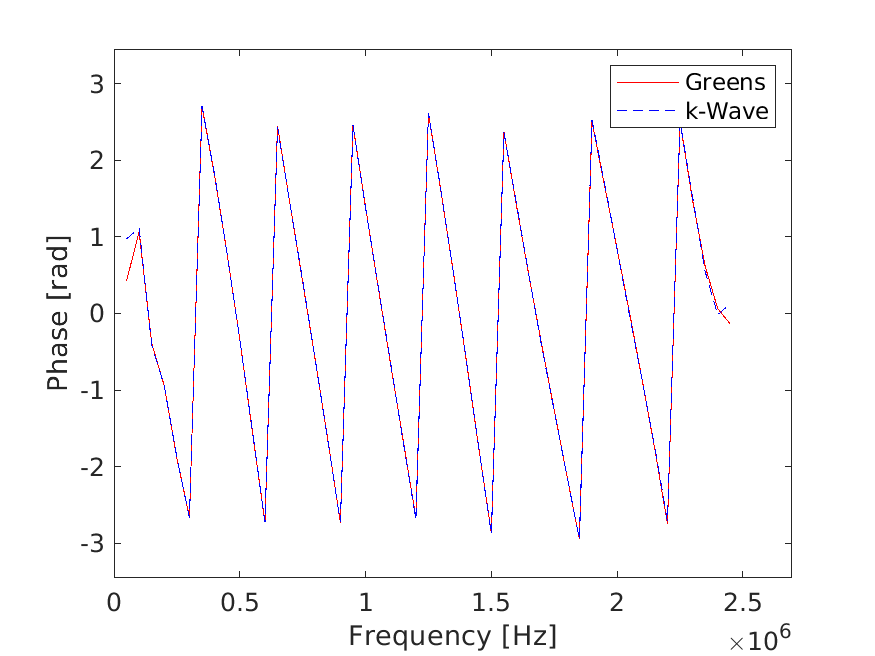}
   \label{fig:signal_water_phase}  }
	\caption{The frequency domain (a) amplitude and (b) phase of the pressure signal measured by receiver 100 after being produced by emitter 1 and propagating through water. (Note that the amplitude has the same normalisation as the excitation signal shown in Fig.\ \ref{fig:excitation-pulse}.)}\label{fig:signal_water}
\end{figure}

\begin{figure} 
   	 \subfigure[]{\includegraphics[width=0.30\textwidth]{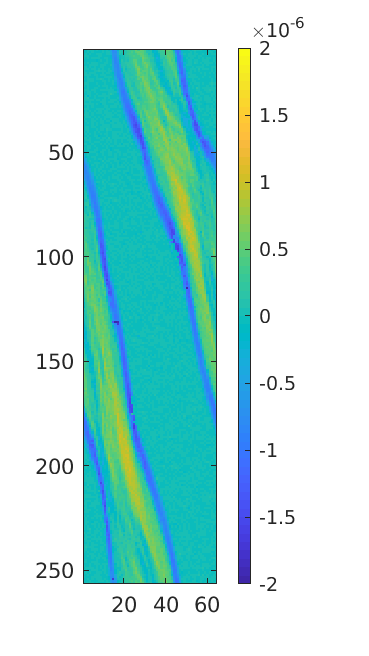}
     \label{fig:tofs-kwave}} 
 	 \subfigure[]{\includegraphics[width=0.30\textwidth]{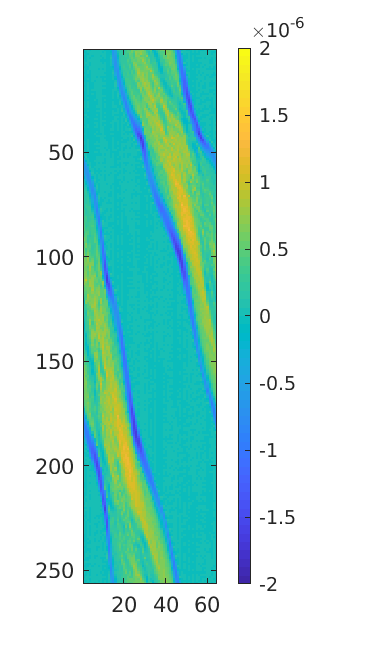}
    \label{fig:tofs-greens}}
	\caption{(a) the discrepancy of the time-of-flight of two sets of time traces measured from the pressure field propagating through the breast inside water and only water.  The time-of-flights were computed using a first-arrival picking algorithm \cite{Li3}, and are approximate. (b) computed $(\phi -\phi_0)/\omega$, which is the same for all frequencies for $y=2$. Here, $\phi_0$ and $\phi$ are accumulated (unwrapped) phases in water and phantom.}\label{fig:tofs}
\end{figure}

\begin{figure} 
   	 \subfigure[]{\includegraphics[width=0.80\textwidth]{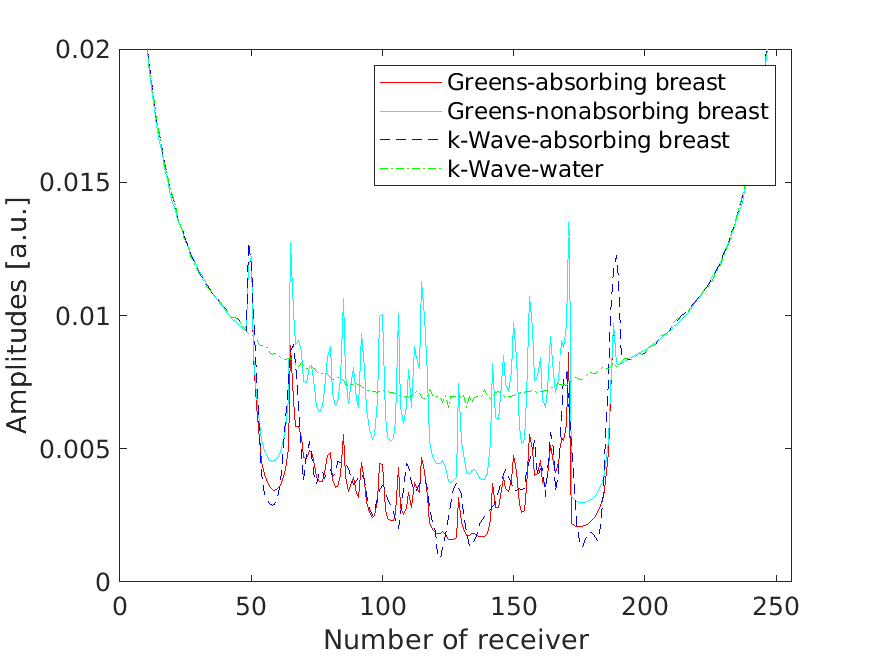}
     \label{fig:amplitudes}} 
	\caption{The amplitudes of the pressure field induced by an excitation  of emitter 1 and measured on the receivers. Approximate Green's function method with absorbing breast phantom (red) and non-absorbing (light blue), and k-wave with absorbing breast phantom (blue) and through water (green). The exponent power $y$ was set 2 for this comparison.}\label{fig:amplitudes}
\end{figure}

\subsection{Approximate Green's forward model} \label{sec:results-greens}
In this section, our approximate Green's function solution to Szabo's the wave equation \eqref{eq:second-order-lossy}, proposed in Sec.\ \ref{sec:greens}, is compared to a full-wave solution using k-wave \cite{kWave,Treeby,Treeby2}. Note that k-wave replaces the fractional time derivative in 
\eqref{eq:second-order-lossy} with two fractional Laplacian operators to improve the computational efficiency \cite{Treeby2}.

\subsubsection{Homogeneous medium}  \label{sec:num-greens-hom}
The pressure $p(\omega, \bx; \bx_e) $ was computed using \eqref{eq:pressure-greens-source}, in which the pressure source was simulated as in Fig.\ \ref{fig:excitation-pulse}, and the Green's function was calculated analytically using the formula \eqref{eq:greens-homogeneous}.
The computation was done at 50 equidistant discretised frequencies $\omega \in  0,...,2 \pi f_{max}  $, where $f_{max} $ is the maximum frequency supported by the computational grid for the k-wave simulation. 
Figs.\ \ref{fig:signal_water_amplitude} and \ref{fig:signal_water_phase} show the amplitude and phase components of the pressure signal measured by receiver 100 after propagation from emitter 1 through water (no breast phantom). The relative discrepancy between the k-wave and the ray-based Green's approach for the displayed signal was $0.75 \% $. For the pressure field produced by emitter 1, the mean relative discrepancy between the k-wave and Green's approach among the signals recorded by all the 256 receivers was $0.77 \% $. The major portion of the discrepancy corresponds to the interpolation from the emission points to the grid points, and the interpolation from the grid points to the reception points.

\subsubsection{Heterogeneous medium} \label{sec:num-greens-het}
This section compares the forward model using our approximate Green's function approach and k-wave for the case of heterogeneous media. The pressure field $p(\omega, \bx; \bx_e)$ induced by each emitter, and propagated through the digital breast phantom, was recorded by all 256 receivers. The simulation for k-wave was performed as explained above.
The pressure field was also approximated using \eqref{eq:pressure-greens-source}, in which the Green's function $g(\omega, \bx; \bx_e) $ is approximated using \eqref{eq:greens-heterogeneous-sampled} along the linked ray, as described in Sec. \ref{sec:Discretised-approximation}.

We first compare the k-wave solver with the approximate Green's function approach for simulating the pressure field on the receivers. For comparison, the k-wave was implemented for simulating the pressure field on the receivers, as described above, and the forward model using the Green's function approach was implemented on the same computational grid and with the same parameters as used for image reconstruction. In the image reconstruction below, the inverse crime was avoided by implementing the Green's function forward model on a grid different from the grid used for data simulation. This grid has a size $200 \times 200 $ with a grid spacing of $1 $mm.
For implementing the Green's approach, the wave number field has been smoothed by an averaging window of size $7$ grid spacings. Only in this section, for making the comparison fair, the wave number field of the digital breast phantom was smoothed by an equivalent averaging window for the grid for data simulation as well. (Note that for simulating the data used for image reconstruction, the k-wave was applied on the nonsmooth digital breast phantom.) The reason was that the approximate Green's forward model neglects the scattering effects, but the primary scattering effects will be included in the image reconstruction through formation of the linearised subproblems \eqref{eq:linearised-subproblem2} from low-to-high frequencies.

For both comparison of forward models and image reconstruction, the step size $\Delta s$ for ray tracing was chosen the same as grid spacing.
The ray linking for 2D UST is equivalent to solving a one-dimensional optimisation problem \cite{Anderson2}.
The most popular one-dimensional approaches for ray linking in 2D UST are \textit{Regula Falsi} and \textit{Secant} methods. The \textit{Regula Falsi} approach is more time consuming, but is more robust, and more likely gives a stable solution than the \textit{Secant} approach when a good initial guess is not available. On the other hand, \textit{Secant} method is faster and is a good choice when the initial guess is sufficiently close to the true solution. Only in this section, because ray linking was applied on the digital breast phantom for the purpose of comparison with k-wave, a good initial guess is not available.  Therefore, \textit{Regula Falsi} approach was used for ray linking for improving the stability \cite{Anderson2}.

The accumulated phase along the linked rays was computed using \eqref{eq:phase-discretised}. Because the rays are linked between emitters and receivers, the receivers match the last point along the rays, and therefore, the pressure can be approximated on the receivers without any interpolation. 
The changes in the amplitude because of cumulative acoustic absorption and geometrical spreading were computed along the linked rays using \eqref{eq:amplitude-abs-discretised} and \eqref{eq:amplitude-geom-discretised}, respectively. For the latter, the amplitude on the first point along the ray was calculated using the analytical formula \eqref{eq:amplitude-homogenous}, and using an assumption that a neighborhood of each emission point with a radius larger than the ray spacing is acoustically homogeneous. The Jacobian on the rays' points was computed using \eqref{eq:derivative_angle-reference}, which required two auxiliary rays for each linked ray (emitter-receiver pair).

 Fig. \ref{fig:tofs-greens} shows the discrepancy $\big(\phi - \phi_0  \big)/ \omega $ using an exponent power $y=2$ (used only for the comparison case), which makes the plot independent from frequencies. Here, $\phi$ represents the accumulated phase along the linked rays for the digital breast phantom, and was approximated using Eq. \eqref{eq:phase-discretised}, and $\phi_0$ is the phase on the receivers in the only water, and was calculated using Eq. \eqref{eq:phase-homogenous}. For the simulated time traces using k-wave, the phase which can be computed is not accumulated, and therefore, the sinogram $\big(\phi - \phi_0  \big)/ \omega $ will suffer from the discontinuities at $-\pi$ and $\pi$. As a result, an approximate broadband variant of the accumulated (unwrapped) phase of time traces was determined using the time-of-flights computed using a first-arrival picking algorithm often used for transmission imaging for medical applications \cite{Li3}. As described, using $y=2$, the effects of dispersion on the plots are neglected. (Note that for the data simulated for image reconstruction, a more realistic $y=1.4$ was used.) Fig. \ref{fig:tofs-kwave} shows the discrepancy of the time-of-flights in the two sets of time traces measured on the receivers, i.e, propagating in water and the digital breast phantom inside the water. (It is reminded that the measured signals imply the simulated signals using k-wave.) 
 In this figure, the columns (resp. rows) represent the emitters (resp. receivers).

Fig. \ref{fig:amplitudes} compares the the k-wave and the Green's approach for the amplitude of the pressure field induced by emitter 1 and measured on all the receivers at frequency $f= 1$MHz. Note that in contrast to the approximate Green's function approach, which approximates the amplitude along the linked rays, the k-wave simulation also includes the scattering effects, but the smoothing applied on the digital breast phantom (only for comparison purposes), reduces the scattering effects.

The blue (resp. green) line represents the recorded pressure on the receivers, when the pressure field propagating through the breast phantom (resp. only water) was simulated using the k-wave. The Green's approach was then used for approximating the pressure field on the receivers. The red (resp. light blue) line represents the amplitude of the Green's function on the receivers for the absorbing breast phantom (resp. the breast phantom with $\alpha=0$). Note that for the light blue plot, only the geometrical effects have been accounted for in calculating the amplitudes, ie. $A_{abs}=1$. It is also reminded that the amplitudes calculated by the analytic Green's function given in Eq. \eqref{eq:greens-homogeneous} match the amplitudes simulated by k-wave in water (green line), as shown in Fig. \ref{fig:signal_water_amplitude}.

\subsection{Image reconstruction}\label{image_reconstruction}
This section gives further details about the procedure used for image reconstruction, and shows the reconstructed images. The grid for image reconstruction consisted of $200 \times 200$ grid points with position $\big[-9.98, +9.93\big] \times \big[-9.98, +9.93\big] \text{cm}^2$ and a grid spacing of $ \Delta x = 1 $mm along all the Cartesian coordinates. The pressure time series were downsampled by two, providing a sampling rate of $26.4$ MHz ($3.7879 \times 10^{-8} \ s$ time spacing) for the data used for image reconstruction. Because the sound speed is reconstructed on the grid points, the parameters of Green's functions approximated along the linked rays must be interpolated to the grid points. This was done by enforcing a triangulation to the sampled points on the linked rays, and then interpolating the approximated parameters on the rays' points to the grid points using a trilinear interpolation. For computing the Green's function $ g(c; \bx_r; \bx) $, the parameters are reversed along the linked rays. This is equivalent to replacing $e$ and $r$ in formulae given in Sec. \ref{sec:Discretised-approximation}. Note that for computing the geometrical portion of the amplitude along each reversed ray using Eq. \ \eqref{eq:derivative_angle-reference}, two additional auxiliary rays must be computed with initial positions $\bx_r$. The accuracy of the reconstructed images are measured in terms of \textit{Relative Error} 
\begin{align}
    RE \left( c_{{image}}\right) = \frac{\| c_{{image}} - c_{{phantom}} \|}{\| c_0 - c_{{phantom}} \|} \times 100,
\end{align}
where $c_{image}$ denotes the reconstructed sound speed distribution. This parameter is calculated on the grid, and inside the binary mask, used for image reconstruction.

\subsubsection{Initial guess.} 
Because the inverse problem of reconstructing the sound speed image from boundary pressure data is nonlinear \cite{Plessix}, an initial guess using a time-of-flight approach is often used \cite{Huth,Wang,Matthews,Matthews1, Liva,Wiskin1,Wiskin2}. Here, an image reconstruction approach based on the {time-of-flight} (TOF) of the measured pressure data \cite{Li3,Javaherian} was used in order to provide an initial guess for the proposed inversion approach, but a rough initial guess was chosen using only the early iterations of the TOF-based inversion in order to indicate that the success of the proposed ray-based Green's inversion approach is not strongly dependent on the image provided by the TOF-based inversion approach.

The discrepancy of the first-arrival of the measured (simulated) pressure time series for the two data sets simulated for the breast phantom inside water and only water was calculated using a first-arrival picking algorithm \cite{Li3}. The TOF-based inversion approach iteratively minimises the norm of discrepancy of first-arrivals modelled by the ray tracing algorithm \ref{alg:ray_tracing_mixed_step} and those calculated from the measured data sets using the first-arrival picking algorithm \cite{Li3}. The sound speed distribution was initialised from the sound speed in water, and the time-of-flights are iteratively modelled as the integral of the slowness (reciprocal of sound speed) along the linked rays. At each iteration, the ray linking was done separately for all emitter-receiver pairs by an iterative implementation of the ray tracing algorithm \ref{alg:ray_tracing_mixed_step} using the \textit{Secant} method, and the initial angle of each ray for ray linking was chosen the optimal initial angle obtained from ray linking at the last iteration \cite{Javaherian}. Only for ray tracing associated with the TOF-based inversion approach, the dispersion effects were ignored, ie. $\alpha=0$. For ray tracing and construction of the system matrix at each iteration (cf. \cite{Javaherian}), the sound speed updates were smoothed by an averaging window of size 7 grid spacings. The optimal TOF-based image was obtained after 7 iterations, ie. one iteration using straight rays following 6 iterations using bent rays applied on the updates of the sound speed map. The first iteration was done using straight rays, because the initial guess was chosen homogeneous, ie. the sound speed in water. Fig. \ref{fig:tof_trial2} shows the image reconstructed after the first iteration of the TOF-based inversion approach using only straight rays ($RE=85.66 \%$).
Fig. \ref{fig:tof_trial1} shows the TOF-based reconstructed image after one iteration using straight rays and subsequently 3 iterations using the bent rays ($RE=74.83 \%$).
Fig. \ref{fig:tof_optimal} shows the optimal image reconstructed by the TOF-based inversion approach after the 7 iterations ($RE=65.16 \%$).

\subsubsection{Green's function-based inversion approach}
The inversion approach explained in Sec.\ \ref{sec:optimisation-greens} was implemented at 140 equidistant discretised frequencies $ f \in {0.2,..., 1.5}$ MHz. The image reconstruction was performed from low to high frequencies such that each linarised subproblem \eqref{eq:linearised-subproblem} was solved at four consecutive discretised frequencies using Algorithm  \ref{alg:conjugate-gradient}. 
For each linearised subproblem, the Green's functions were approximated along the rays traced in a medium having the last update of the sound speed. The ray linking was performed separately for all emitter-receive pairs using the \textit{Secant} method. For ray linking, a set of constraints was enforced on the initial angles in order to improve the stability, as explained in \cite{Javaherian}. The trajectory of rays are iteratively computed using Algorithm \ref{alg:ray_tracing_mixed_step} until the end point of the rays matches the position of the reception points within a tolerance. For computing the trajectory of rays, an averaging window of size $7$ grid spacings was enforced on the updated wave number maps, but the nonsmoothed maps were used for integration along the rays and approximating the Green's functions using the formulae given in Sec. \ref{sec:Discretised-approximation}. For each emitter-receiver pair, the initial angle of the linked ray for the linearised subproblem $n$ was used as the initial guess for ray linking for the linearised subproblem $n+1$ \cite{Javaherian}. For $n=0$, the initial guess for the initial angles were chosen as those obtained from ray linking at the last iteration of the TOF-based inversion approach. Each linearised subproblem was solved using maximum $l_{max} =10$ cg iterations (cf. Algorithm \ref{alg:conjugate-gradient}). 

\begin{figure} 
    \centering
     \subfigure[]{\includegraphics[width=0.45\textwidth]{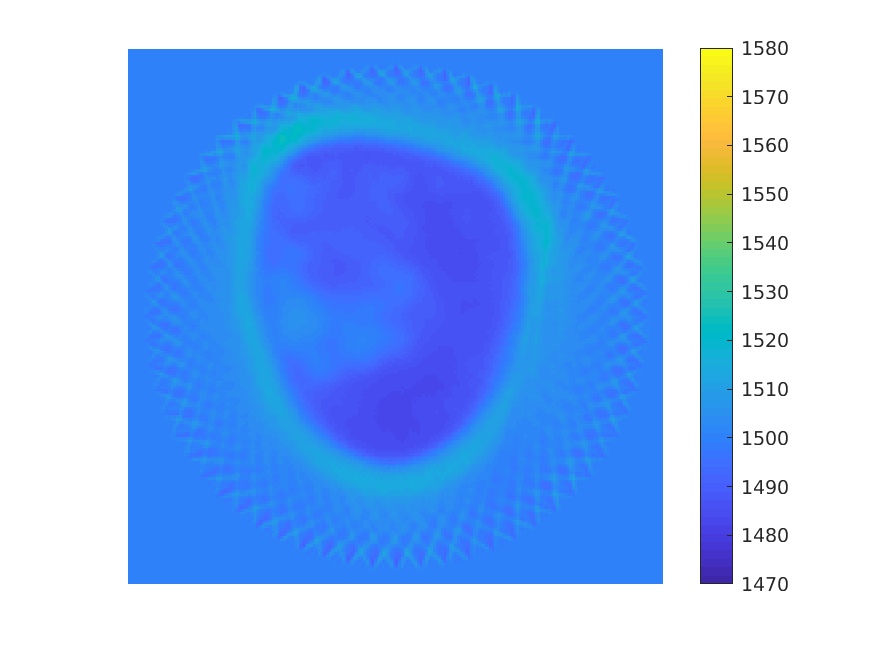}
     \label{fig:tof_trial2}  }
      \subfigure[]{\includegraphics[width=0.45\textwidth]{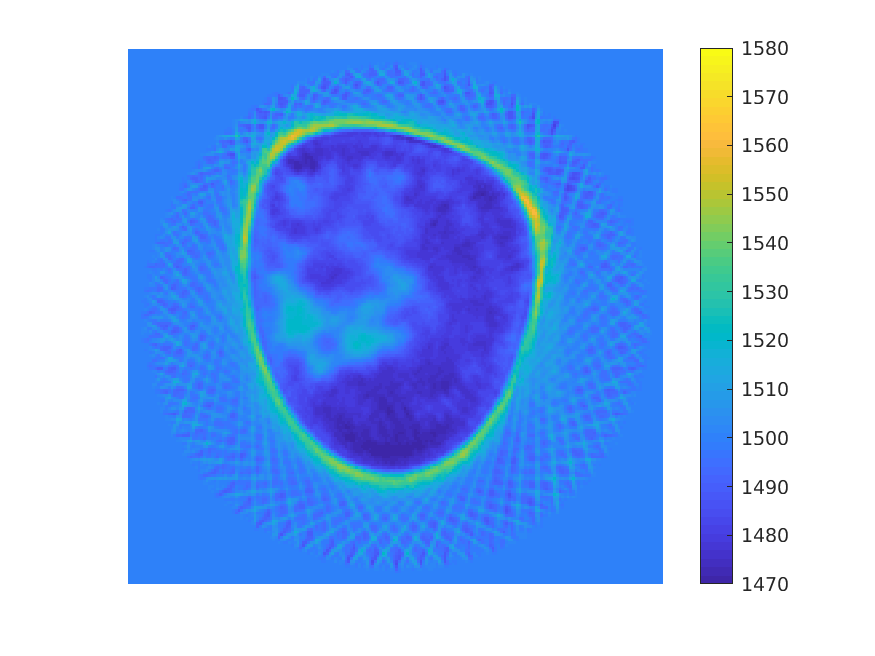}
     \label{fig:tof_trial1}   }
     \subfigure[]{\includegraphics[width=0.45\textwidth]{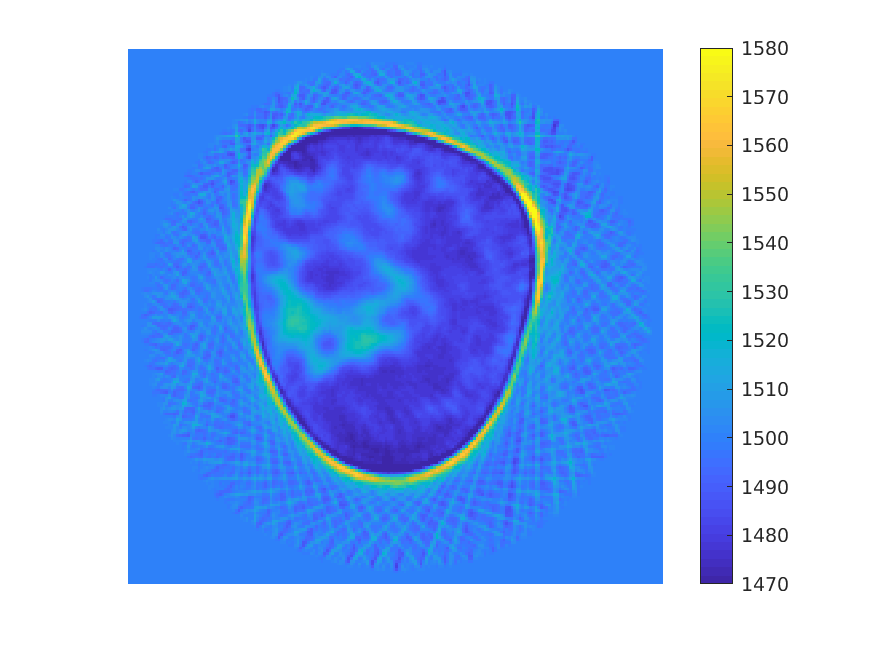}
     \label{fig:tof_optimal}   }
	\caption{TOF-based reconstructed images after: (a) one iteration using straight rays (b) one iteration using straight rays and 3 iterations using bent rays  (c) one iteration using straight rays and 6 iterations using bent rays (the optimal image).}
\end{figure}

\begin{figure} 
    \centering
      \subfigure[]{\includegraphics[width=0.45\textwidth]{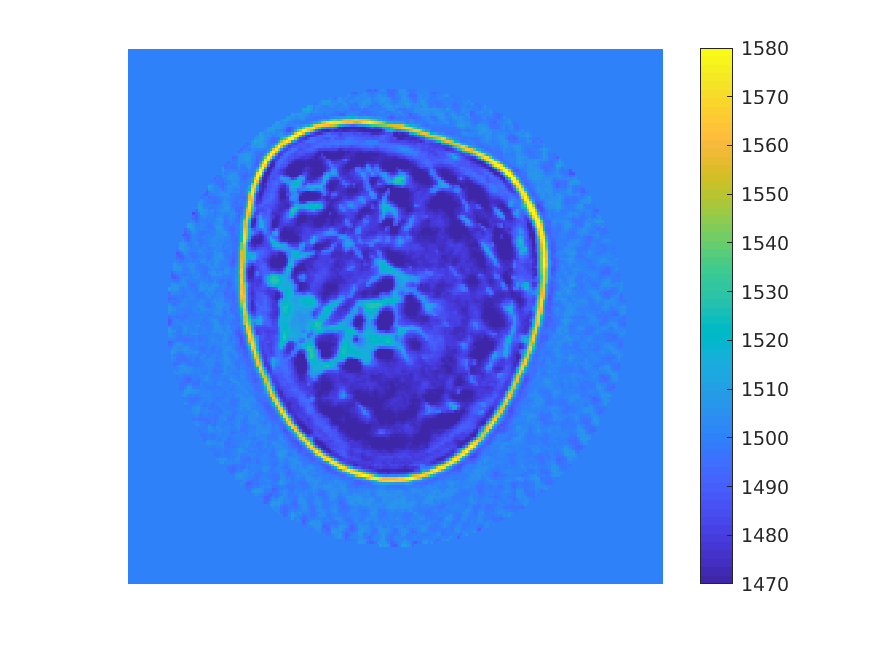}
     \label{fig:greens_trial2}   }
      \subfigure[]{\includegraphics[width=0.45\textwidth]{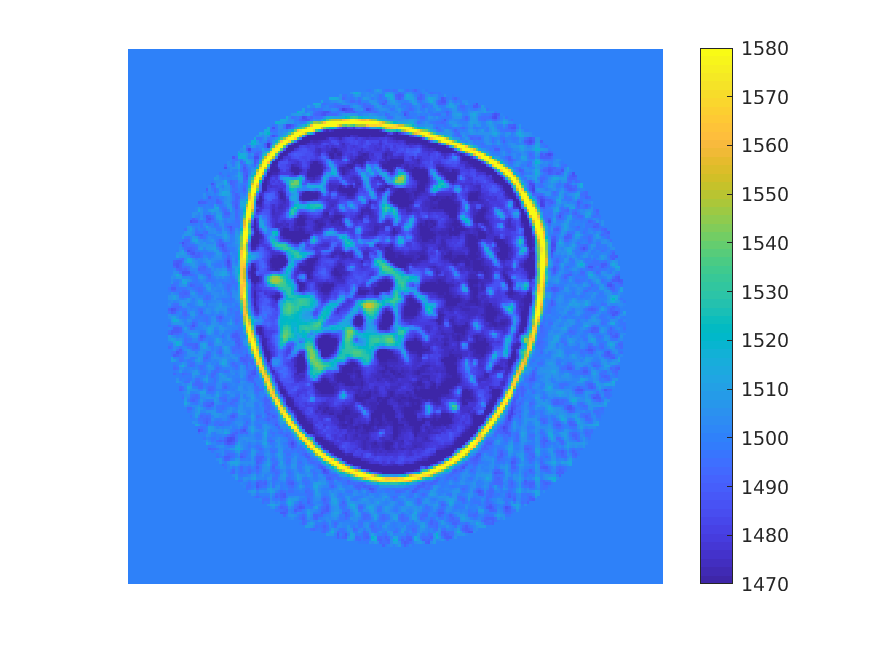}
     \label{fig:greens_trial4}   }
      \subfigure[]{\includegraphics[width=0.45\textwidth]{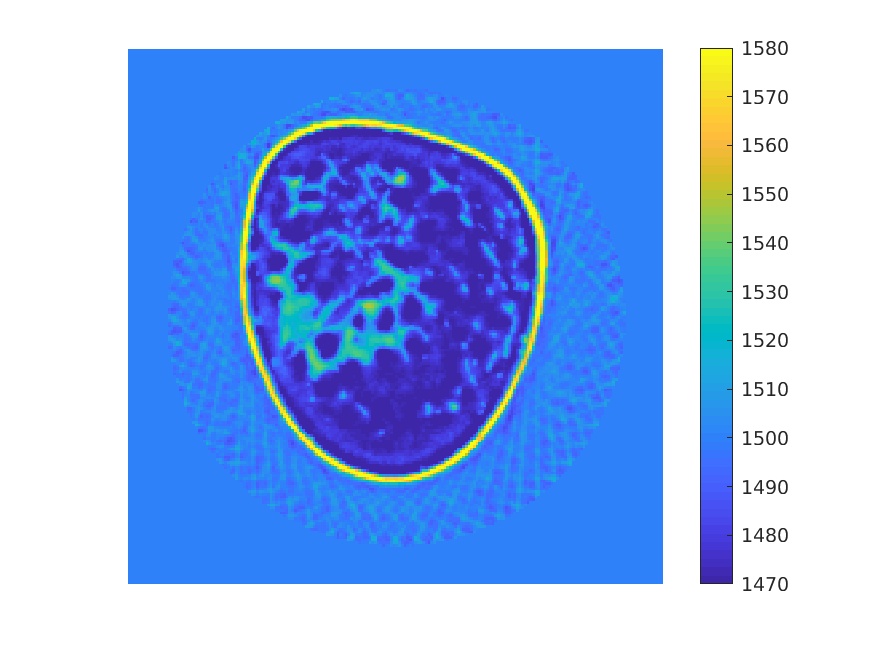}
     \label{fig:greens_trial3}   }
      \subfigure[]{\includegraphics[width=0.45\textwidth]{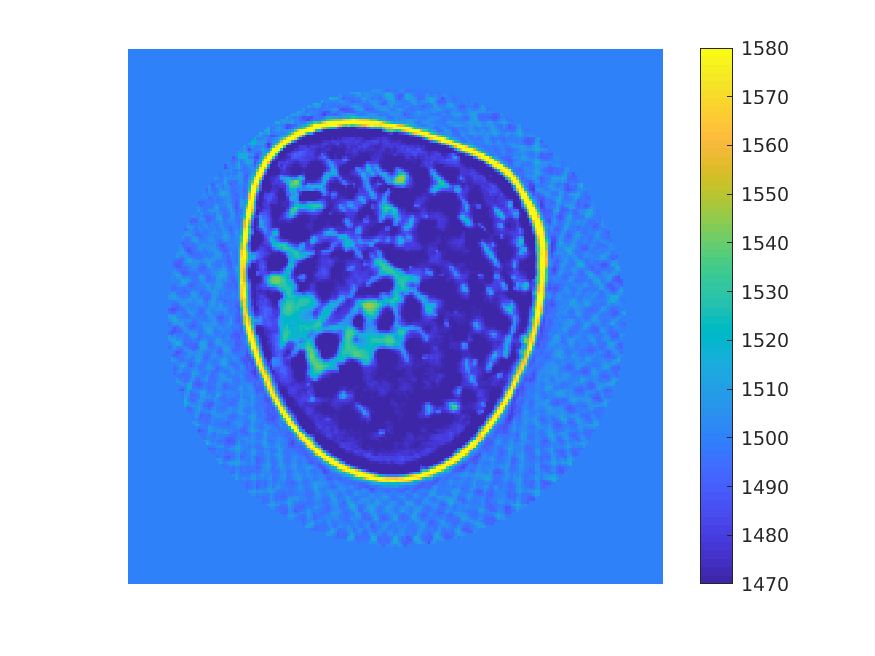}
     \label{fig:greens_trial1}   } 
	\caption{Reconstructed images using the ray-based Green's approach. (a) initial guess far from the ground truth (Fig. \ref{fig:tof_trial2}), exact $\alpha_0$. Initial guess closer to the ground truth, but not the optimal (Fig. \ref{fig:tof_trial1}): (b) $\alpha_0 = 0$, (c) homogeneous $\alpha_0 = 5 \ \text{dB} \text{MHz}^{-y} \text{cm}^{-1}$ inside the phantom, (d) true $\alpha_0$ (Fig. \ref{fig:phantom_absorption_coefficient}).}
\end{figure}

Fig. \ref{fig:greens_trial2} shows the sound speed image reconstructed using the proposed Green's function approach ($RE=45.93\%$), when the initial guess was chosen to be the sound speed estimate shown in Fig. \ref{fig:tof_trial2}, and the true absorption coefficient map shown in Fig. \ref{fig:phantom_absorption_coefficient} was used for image reconstruction.
Figs. \ref{fig:greens_trial4},  \ref{fig:greens_trial3} and \ref{fig:greens_trial1}
show the sound speed images reconstructed using the Green's function approach, when the initial guess was chosen to be the sound speed shown in Fig. \ref{fig:tof_trial1}. Fig. \ref{fig:greens_trial4} shows the sound speed image ($RE=42.51\%$), when the absorption and dispersion are neglected by using $\alpha_0 = 0$. 
Fig. \ref{fig:greens_trial3} shows the reconstructed sound speed image ($RE=38.87\%$), when the absorption coefficient inside the breast phantom was set homogeneous and $\alpha_0 = 5 \ \text{dB} \text{MHz}^{-y} \text{cm}^{-1}$ in order to avoid an inverse crime. Fig. \ref{fig:greens_trial1} shows the reconstructed image ($RE=37.91\%$), when the true absorption coefficient given in Fig \ref{fig:phantom_absorption_coefficient} was used.
As shown in these images, the discrepancy between the reconstructed images using the true $\alpha_0$ and the erroneous homogeneous absorption coefficient map is small, but the reconstructed image for which the absorption and dispersion effects are neglected includes more artefact.
It is worth mentioning that for choosing the initial guess of the Green's inversion approach as a sound speed image reconstructed using the TOF-based approach, a trade-off should be made between the closeness to the optimal solution of the Green's approach and the artefact included in the reconstructed image \cite{Wiskin1, Wiskin2}. 

For purposes of comparison, a reconstruction was performed using only straight rays (not shown). We found that using bent rays resulted in a slightly more accurate sound speed update than using straight rays at low frequencies (200 kHz), but the improvement became greater as the frequency increased.

\section{Discussion}
\label{sec:discussion}

In this section we highlight some similarities with other approaches used to tackle UST, and briefly discuss the advantages and disadvantages of this approach compared to other commonly-used approaches. There is a very extensive literature on UST dating back decades, and, for reasons of space, we cannot reference all possible relevant works here, but we hope those works referenced are representative and give the reader a route into the considerable literature.

\subsection*{Inversions using full-wave models}
The most general approach for reconstructing the sound speed from acoustic pressure data is based on minimising the difference (in some sense) between the measurements and a model of the acoustic propagation. This model-based minimisation framework is widely used for tackling inverse problems, and has also been widely used for UST. It is also the approach taken in this paper. Arguably, the choice of model that will achieve the most accurate results is the model that describes the physics of acoustic propagation most accurately. For this reason, there has been great interest recently in using so-called full-wave models that explicitly model the acoustic wave equation for heterogenous media \cite{Wang,Matthews,Matthews1,Goncharsky3, Liva,Bachmann,Lucka}. As expected, these approaches have been shown to provide accurate and high resolution images, but a significant hurdle to the practical applicability of these is that the full-wave solvers are compute intensive. Indeed, grid-based solvers become increasingly memory-hungry the higher the frequency. (Another potential disadvantage is when using a rotating measurement system \cite{Ruiter,Gemmeke}, as the number of forward solves grows linearly with the number of rotations \cite{Javaherian}.)
It is also important to consider that it is never possible to include every conceivable physical effect in any forward model, and that even full-wave models are based on approximations. For example, it is not necessary to model backscattered waves when measurements are only made in front of the emitter. The challenge is balancing the extent to which the physics is modelled accurately with the cost of computation, such that good images can be produced in a reasonable time. For this reason, approximate models have been proposed for use within model-based inversions
\cite{Borup,Wiskin1,Wiskin2}. The approach proposed here falls into this category. Our approximate forward model is based on a ray-approximation to the Green's function. This has two significant advantages over full-wave models: it is much more computationally efficient than grid-based full-wave solvers, and the cost does not scale with frequency. Also, it is trivial to incorporate any arbitrary transducer element directionality in a ray-based scheme (for transmit and receive), which is not easy to do in full-wave models and can be important in practice. Now, it is well known that ray theory is a very good approximation when the frequency is sufficiently high that it sees the medium as smoothly varying. Here, we showed that this approximation, when embedded in a second-order minimisation scheme, can recover very accurate images for the level of contrast found in soft tissue.

\subsection*{Inversions based on time-of-flight}
Ray-based inversion schemes for UST have typically been used when the data is given in the form of a set of times-of-flight between the emitters and receivers. This computationally efficient approach has been widely used in 2D \cite{Anderson1,Anderson2,Anderson3,Duric,Li1} and more recently in 3D \cite{Javaherian}. In our approach here, we retain the computational efficiency of ray-tracing, but dispense with the drawbacks of time-of-flight picking by using the whole set of acoustic time series in the inversion. Time-of-flight picking, by reducing each time series to a single number, removes a great deal of information in the data, effectively reducing the signal-to-noise ratio. 

It might be objected that a ray-based forward model, such as that used here and included in the objective function $\mathcal{F}$, does not account for scattering, only refraction, geometric spreading and absorption. While that is true, a key point of our approach is that the Fre\'chet derivative, included in both the gradient and Hessian, does take into account the primary scattered waves, so these are included during the image reconstruction (over an increasingly good estimate of the sound speed at each iteration). This explained in the context of the Born approximation below.

Note that we use ray-tracing to compute the rays, rather than grid-based eikonal solvers such as fast marching, as the latter suffer from similar computational challenges to grid-based full-wave models, especially in a rotational setting and with a sparse distribution of receivers.

\subsection*{Inversions exploiting the Born-approximation}

A common approach to incorporating first-order scattered waves in simulations and inversions is to turn to the Born approximation, and in particular the ray-Born approximation \cite{Coates,Mojabi,Sarajaervi}. The Born approximation suggests an efficient way for approximating the wavefield using an assumption that the sound speed model can be decomposed into a slowly varying background and a rapidly varying scatterer \cite{Moser}. It has been studied extensively from both theoretical \cite{Hudson1,Hudson2} and numerical points of view \cite{Cerveny2001,Sarajaervi}. 

Born forward modelling assumes that the wavefield measured on the receivers consists of two components: the free-space Green's function solution in the homogeneous or slowly varying background medium, and a scattering integral, which accounts for the scattered waves up to the order of a scattering series used for describing the relation between the perturbation in the sound speed and the induced perturbation in the wavefield. This scattering series is referred to as \textit{Born series} \cite{Moser}. The first-order approximation of the Born series linearly relates a perturbation to the sound speed to the primarily scattered waves.  

The Born approximation was used for including the diffraction of plane waves in the transmission computed tomography \cite{Devaney,Devaney2}. In the context of diffraction tomography, the Fourier diffraction theorem, which relates the Fourier transform of the scattered wavefield to the Fourier transform of the medium's sound speed, can be derived using a first-order Born approximation \cite{Kak}. In one clinical application, it was shown that transmission diffraction tomography relying on the Born approximation provided higher resolution images than time-of-flight-based transmission ultrasound tomography \cite{Simonetti}. However, the Born approximation is only valid approximation under certain conditions that do not hold for transmission imaging of soft tissues like the breast, ie. the variations of the sound speed about the background medium must be very small, and the object's size must be small sufficiently that the distortion in phase be smaller than $\pi$ \cite{Kak,Huth}. To ameliorate these problems, a time-of-flight-based image was used to correct for the aberration that causes diffraction tomography to fail under the Born approximation \cite{Huth}.

In the context of seismic imaging, the class of ray-based migration/inversion approaches, which provide an efficient way for imaging of the scatterers and reflectors, rely on the Born approximation \cite{Coates,Thierry,Lambare1,Lambare2}. Additionally, the frequency-space inversion approaches based on the Born approximation have been widely used in seismic \cite{Pratt} and medical applications \cite{Wiskin1,Wiskin2}. 

Here, the Born approximation was used for an iterative linearisation of an objective (misfit) function. The objective function is defined in the frequency domain, and is minimised from low-to high frequencies. The problem associated with large phase distortion or the large sound speed variations, which limits the application of the Born approximation to soft tissue, was addressed by embedding it in a second-order optimisation framework. The nonlinear objective function was iteratively minimised through a series of linearisations using the Born approximation. This approximation of the objective function gives a linear relation between the perturbation in the sound speed and the induced perturbed pressure field. For each linearsation, the background medium is chosen the last update of the sound speed from the previous linearisation (frequency range), and ray theory was used to approximate the Green's functions. The objective function was iteratively linearised and minimised, starting from low to high frequencies, and the singly scattering features are therefore included in the background medium from low-to-high frequencies. Because of the gradual inclusion of the scatterers in the iteratively-found background medium, the problem of the multivaluedness of the ray field in the background medium may occur \cite{Moser}. To tackle this, either more advanced asymptotic modeling techniques such as Gaussian beams should be used \cite{Cerveny1983} at the cost of more computational cost, or the background should be smoothed \cite{Moser}. Here, the latter approach was used such that the rays' trajectory was computed on the sound speed updates smoothed by an averaging window, which was kept the same for all frequencies. (Note that the phase and amplitudes along the rays were calculated using the nonsmoothed version of the sound speed updates using the formulae in Sec.  \ref{sec:Discretised-approximation}). 

Our numerical experience shows that even with applying the averaging window on the sound speed (or wavenumber) field, the ray linking is the key step in approximating accurate pressure time series at the receivers.  Avoiding ray  linking by tracing  the  rays  using  equidistant  initial  angles,  and  then  interpolating to  the  reception  points  will  introduce errors in the acoustic pressure on some of the receivers. This  is  because  a  small  change  in  the  initial  angle  may  lead  to a large change in the rays’ trajectory, and therefore the assumption of differentiability required for the interpolation does not hold.

\subsection*{Absorption} 
 Absorption and dispersion effects were included in the simulated data and in the inversion approach. This study was limited to an image reconstruction of the sound speed given the absorption coefficient map, but it was shown that the proposed method can reconstruct an accurate sound speed map for soft tissues like the breast, when only partial information about the absorption coefficient map is available, ie. the absorption coefficient map across the breast phantom is assumed homogeneous. The quality of the reconstructed image degrades slightly when the absorption and dispersion effects are fully neglected. It should also be noted that the inversion framework presented here will allow for the extension of this approach to recovering both the absorption and the sound speed simultaneously, because both affect the complex scattering potential $\Upsilon \, \delta k$.

\subsection*{Practical applicability}
 Although here, the effectiveness of the proposed inversion approach was demonstrated in a simulation scenario, the inverse crime was avoided by using two inherently different approaches for data simulation and for the forward model in the image reconstruction. Also, the excitation pulse used in this study was the output of a physical transducer, and the number of excitations (emitters) was less than is often used in practice for 2D UST \cite{Duric}, suggesting it will be possible to achieve similar results with experimental data. 
 
 For nonlinear inverse problems, regularisation is often used for reducing the ill-posedness and improving the stability \cite{Matthews,Liva}. Here, although explicit regularisation was avoided, the cg iterations $l_{max}$ implicitly enforce a regularisation on the solution. Furthermore, the use of increasing frequencies (as is common in full-wave inversion approaches and is akin to a multigrid method) helped ensure the solution did not become stuck in a local minimum.
 
 The proposed approach was proposed and demonstrated for a 2D scenario, but the acoustic waves actually travel in 3D medium. An extension of the proposed approach to 2.5D, ie., imaging within a slab-like volume containing a target slice along the detection ring \cite{Li2}, or full-3D \cite{Ruiter,Gemmeke}, is straightforward, and will be studied in future work. A key challenge in extending this approach to 3D is 3D ray-linking; we recently proposed an efficient and robust approach to overcome this \cite{Javaherian}.

\section{Summary}
An efficient UST image reconstruction algorithm is described and demonstrated for recovering the sound speed distribution from acoustic time series measurements made in soft tissue. The approach is based on a second-order iterative minimisation of the difference between the measurements and a model based on a ray-approximation to the heterogeneous Green's function.

\section*{Acknowledgements}
This work was funded by the European Union’s Horizon 2020 Research and Innovation program H2020 ICT 2016-2017 under Grant agreement No. 732411, which is an initiative of the Photonics Public Private Partnership. The authors would like to thank several people for helpful discussions around the topic of this paper. In particular, Marta Batcke and Francesc Rul$\cdot$lan on rays, Felix Lucka on UST, and Bradley Treeby on k-Wave.

\section*{Appendix: Derivation of Hessian matrix}

Here, further details will be given about the Hessian matrix (cf. Eq. \ \eqref{eq:hessian}). Considering \eqref{eq:linearised-subproblem}, the Hessian matrix satisfies
\begin{align}
 H^{(n)} (\bx, \bx') = 
 \frac{  \partial   }{\partial  c^{(n)} (\bx) } 
\Big[ \frac{\partial \mathcal{F} ( c^{(n)} )}{\partial \big( c^{(n)} (\bx') \big)} \Big]=
\frac{  \partial   }{\partial \big( c^{(n)} (\bx) \big)} \nabla \mathcal{F}^{(n)} (\bx')
 \end{align}
Plugging Eq. \eqref{eq:gradient-greens2} into the right-hand-side of the above equation gives two terms \cite{Margrave}
 \begin{align}
 H^{(n)} (\bx, \bx') = H_1^{(n)} (\bx, \bx') + H_2^{(n)} (\bx, \bx'),
 \end{align}
where the first term is in the form
\begin{align}  \label{eq:hessian_first_term}
\begin{split}
    & H_1^{(n)} (\bx, \bx')  =\\
     &\sum_{e,r}  \int  \text{Re} \left\{     \frac{  \partial   }{\partial \big( c^{(n)} (\bx) \big)} \Big[  \Upsilon_{c}^{(n)} (\bx')     g(c^{(n)}; \bx_r; \bx')   \big[  g(c^{(n)}; \bx'; \bx_e)  s(\bx_e)  \big]    \Big]
         \  \delta p_{res}^* (c^{(n)}; r, e)  \right\} d \omega.
\end{split}
\end{align}
Here,
\begin{align}
\begin{split}
   &\frac{  \partial   }{\partial \big( c^{(n)} (\bx) \big)} \Big[  \Upsilon_{c}^{(n)} (\bx')     g(c^{(n)}; \bx_r; \bx')   \big[  g(c^{(n)}; \bx'; \bx_e)  s(\bx_e)  \big]    \Big]\\
   &= \Big( \frac{  \partial   }{\partial \big( c^{(n)} (\bx) \big)} \Upsilon_{c}^{(n)}(\bx')  \Big)  \  g(c^{(n)}; \bx_r; \bx')   \big[  g(c^{(n)}; \bx'; \bx_e)  s(\bx_e)  \big]\\
   &+\Upsilon_{c}^{(n)}(\bx')  \  \frac{  \partial   }{\partial \big( c^{(n)} (\bx) \big)} \Big( g(c^{(n)}; \bx_r; \bx')   \big[  g(c^{(n)}; \bx'; \bx_e)  s(\bx_e) \big] \Big) ,  
\end{split}
\end{align}
where in the last line in the above equation, 
\begin{align}
\begin{split}
&\frac{  \partial   }{\partial \big( c^{(n)} (\bx) \big)} \Big( g(c^{(n)}; \bx_r; \bx')   \big[  g(c^{(n)}; \bx'; \bx_e)  s(\bx_e) \big] \Big) \\
&= g(c^{(n)},\bx_r;\bx) \Upsilon_c^{(n)}(\bx)  g(c^{(n)},\bx; \bx') \big[  g(c^{(n)}; \bx'; \bx_e)  s(\bx_e) \big]  \\
&+g(c^{(n)},\bx_r;\bx')   g(c^{(n)},\bx'; \bx) \Upsilon_c^{(n)}(\bx) g(c^{(n)},\bx; \bx_e) s(\bx_e),
\end{split}
\end{align}
where the derivatives have been calculated using Eq. \eqref{eq:Frechet-derivative}.

\noindent
Also, the second term $ H_2^{(n)} $ satisfies \label{eq:hessian_second_term}
\begin{align}
\begin{split}
    & H_2^{(n)} (\bx, \bx')  =\\
     &\sum_{e,r}  \int  \text{Re} \left\{  \Big[    \Upsilon_{c}^{(n)} (\bx)     g(c^{(n)}; \bx_r; \bx)  \big[  g(c^{(n)}; \bx; \bx_e)  s(\bx_e)  \big]    \Big]
      \times     \ \frac{  \partial p^*_{(c^{(n)}; r, e)}   }{\partial \big( c^{(n)} (\bx') \big)}    \right\} d \omega,
\end{split}
\end{align}
where we have used $   \partial /\partial \big( c^{(n)} (\bx') \big) \Big(  \delta p_{res} (c^{(n)}; r, e) \Big) = \partial /\partial \big( c^{(n)} (\bx') \big) \Big(  p_{ (c^{(n)}; r, e)} \Big) $.

Considering that the terms in $H_1^{(n)} $ is negligible compared to the second term $H_2^{(n)} $, and also using the reciprocity of Green's function $g(c^{(n)}; \bx_r; \bx)$, gives a very good approximation to the action of the Hessian on the sound speed perturbation, as given in \eqref{eq:hessian}.

\end{document}